\documentclass{article}
\usepackage{graphicx}
\usepackage{amsmath}
\usepackage{amsfonts}
\usepackage{amssymb}
\newtheorem{theorem}{Theorem}

\newtheorem{corollary}[theorem]{Corollary}

\newtheorem{definition}[theorem]{Definition}

\newtheorem{lemma}[theorem]{Lemma}

\newtheorem{proposition}[theorem]{Proposition}
\newtheorem{remark}[theorem]{Remark}

\newenvironment{proof}[1][Proof]{\textbf{#1.} }{\ \rule{0.5em}{0.5em}}

\begin{document}

\title{Instantons and Branes in \\Manifolds with Vector Cross Product}
\author{Jae-Hyouk Lee and Naichung Conan Leung}
\maketitle
\begin{abstract}
In this paper we study the geometry of manifolds with vector cross product and
its complexification. First we develop the theory of instantons and branes and
study their deformations. For example they are (i) holomorphic curves and
Lagrangian submanifolds in symplectic manifolds and (ii) associative
submanifolds and coassociative submanifolds in $G_{2}$-manifolds.

Second we classify K\"{a}hler manifolds with the complex analog of vector
cross product, namely they are Calabi-Yau manifolds and hyperk\"{a}hler
manifolds. Furthermore we study instantons, Neumann branes and Dirichlet
branes on these manifolds. For example they are special Lagrangian
submanifolds with phase angle zero, complex hypersurfaces and special
Lagrangian submanifolds with phase angle $\pi/2$ in Calabi-Yau manifolds.

Third we prove that, given any Calabi-Yau manifold, its isotropic knot space
admits a natural holomorphic symplectic structure. We also relate the
Calabi-Yau geometry of the manifold to the holomorphic symplectic geometry of
its isotropic knot space.

\pagebreak 
\end{abstract}
\tableofcontents

\pagebreak 

\section{Introduction}

The vector product, or the cross product, in $\mathbb{R}^{3}$ was generalized
by Gray (\cite{B-G},\cite{G}) to the product of any number of tangent vectors,
called the \textit{vector cross product} (abbrev. VCP). The list of Riemannian
manifolds with VCP structures on their tangent bundles include symplectic (or
K\"{a}hler) manifolds, $G_{2}$-manifolds and $Spin\left(  7\right)
$-manifolds. We introduce the complex analog of VCP in definition \ref{Complex
VCP}, called the \textit{complex vector cross product} (abbrev. $\mathbb{C}%
$-VCP). We show that there are only two classes of manifolds with $\mathbb{C}$-VCP.

\begin{theorem}
If $M$ is a closed K\"{a}hler manifold with a $\mathbb{C}$-VCP, then $M$ must
be either (i) a Calabi-Yau manifold, or (ii) a hyperk\"{a}hler manifold.
\end{theorem}

We study the geometry of \textit{instantons} which are submanifolds in $M$
preserved by the VCP. Instantons are always absolute minimal submanifolds in
$M$. When an instanton is not a closed submanifold in $M$, we require its
boundary to lie inside a \textit{brane} in order to have a Fredholm theory for
the free boundary value problem. For example, when $M$ is a symplectic
manifold, then instantons and branes are holomorphic curves and Lagrangian
submanifolds in $M$ respectively. These geometric objects play important roles
in understanding the symplectic geometry of $M$.%

\[%
\begin{tabular}
[c]{|l||l|l|l|l|}\hline
$%
\begin{array}
[c]{l}%
\,\text{Manifolds}\\
\,\text{w/ VCP}%
\end{array}
$ & $%
\begin{array}
[c]{l}%
\text{Symplectic\thinspace}\\
\text{\thinspace Manifolds}%
\end{array}
$ & $%
\begin{array}
[c]{l}%
G_{2}\text{-\thinspace}\\
\text{\thinspace Manifolds}%
\end{array}
$ & $%
\begin{array}
[c]{l}%
Spin\left(  7\right)  \text{-\thinspace}\\
\text{\thinspace Manifolds}%
\end{array}
$ & $%
\begin{array}
[c]{l}%
\text{Oriented\thinspace}\\
\text{\thinspace Manifolds}%
\end{array}
$\\\hline
Instantons & $%
\begin{array}
[c]{l}%
\text{Holomorphic\thinspace}\\
\text{\thinspace curves}%
\end{array}
$ & \multicolumn{1}{|l|}{$%
\begin{array}
[c]{l}%
\text{Associative\thinspace}\\
\text{\thinspace submanifolds}%
\end{array}
$} & $%
\begin{array}
[c]{l}%
\text{Cayley}\\
\text{\thinspace submanifolds}%
\end{array}
$ & $%
\begin{array}
[c]{l}%
\text{Open}\\
\text{\thinspace submanifolds}%
\end{array}
$\\\hline
Branes & $%
\begin{array}
[c]{l}%
\text{Lagrangian}\\
\text{\thinspace submanifolds}%
\end{array}
$ & \multicolumn{1}{|l|}{$%
\begin{array}
[c]{l}%
\text{Coassociative\thinspace}\\
\text{\thinspace submanifolds}%
\end{array}
$} & \multicolumn{1}{|c|}{N/A} & Hypersurfaces\\\hline
\end{tabular}
\]

Instantons $A$ in $M$ can be characterized by the condition
\[
\tau|_{A}=0
\]
with $\tau\in\Omega^{r+1}\left(  M,\mathbf{g}_{M}^{\bot}\right)  $. This is
useful in giving a uniform description of the deformation theory of instantons
(see section \ref{1Sec Deform Inst} for details.). For branes $C$ in $M$, its
normal bundle is naturally identified with the space $\Lambda_{VCP}^{r}%
T_{C}^{\ast}$ of VCP forms of degree $r$ on $C$. Infinitesimal deformations of
branes are parametrized by such differential forms on $C$ which are closed.
Furthermore they are always unobstructed, namely the moduli spaces of branes
are always smooth. (See section \ref{1Sec Deform Brane} for details).

\bigskip

For the \textit{complex} analog, there are two type of branes corresponding to
the Dirichlet type and Neumann type boundary value problems for instantons. We
call them \textit{D-branes} and \textit{N-branes} respectively.\label{Table
C-VCP}
\[%
\begin{tabular}
[c]{|l||l|l|}\hline
$%
\begin{array}
[c]{l}%
\,\text{Manifolds}\\
\,\text{w/ }\mathbb{C}\text{-VCP}%
\end{array}
$ & Calabi-Yau manifolds & hyperk\"{a}hler manifolds\\\hline
Instantons$%
\begin{array}
[c]{l}%
\,\\
\,
\end{array}
$ & Special Lagrangian$_{\theta=0}$ & $I$-holomorphic curves\\\hline
N-Branes$%
\begin{array}
[c]{l}%
\,\\
\,
\end{array}
$ & Complex Hypersurfaces & $J$-complex Lagrangians\\\hline
D-Branes$%
\begin{array}
[c]{l}%
\,\\
\,
\end{array}
$ & Special Lagrangian$_{\theta=-\pi/2}$ & $K$-complex Lagrangians\\\hline
\end{tabular}
\]

\bigskip

Using the Riemannian metric on $M$, any closed VCP determines a differential
form $\phi$ as follows:
\[
\phi\left(  v_{1},...,v_{r},v_{r+1}\right)  =\left\langle v_{1}\times...\times
v_{r},v_{r+1}\right\rangle \text{.}%
\]
By transgressing $\phi$, we obtain a two form on the (multi-dimensional) knot
space of $M$:
\[
\mathcal{K}_{\Sigma}M=Map\left(  \Sigma,M\right)  _{emb}/Diff\left(
\Sigma\right)  \text{,}%
\]
where $\Sigma$ is any smooth manifold of dimension $r-1$. This gives a
symplectic structure on $\mathcal{K}_{\Sigma}M$ which is compatible with the
$L^{2}$-metric, namely an almost K\"{a}hler structure. For instance, when $M$
is a three manifold, $\mathcal{K}_{\Sigma}M$ is the space of knots in $M$. In
this case, Brylinski \cite{Bry} showed that $\mathcal{K}_{\Sigma}M$ is indeed
a K\"{a}hler manifold. He also studied relationship between the geometry of
$\mathcal{K}_{\Sigma}M$ and the geometric quantization of Chern-Simons theory.

In general, VCP geometry on $M$ can be interpreted as the symplectic geometry
on $\mathcal{K}_{\Sigma}M$. For example, for a disk $D$ in $Map\left(
\Sigma,M\right)  _{emb}$ which is horizontal in the principal bundle
\[
Diff\left(  \Sigma\right)  \rightarrow Map\left(  \Sigma,M\right)
_{emb}\overset{\pi}{\rightarrow}\mathcal{K}_{\Sigma}M,
\]
determines a map from $D\tilde{\times}\Sigma$ to $M$ and a disk $\hat{D}$
$=:\pi\left(  D\right)  $ in $\mathcal{K}_{\Sigma}M$. If $\hat{D}$ is a
holomorphic disk in $\mathcal{K}_{\Sigma}M$, $D\tilde{\times}\Sigma$ gives an
instanton in $M$. It has an extra property that the induced metric on
$D\tilde{\times}\Sigma$ gives a Riemannian submersion from $D\tilde{\times
}\Sigma$ to $D$. We call such an instanton a \textit{tight instanton.} In
section \ref{Symplectic on Loop} we prove the following theorem.

\begin{theorem}
Suppose $M$ is a Riemannian manifold with a closed differential form $\phi$.
Then we have

(1) $\phi$ is a VCP form on $M$ if and only if its transgression defines an
almost K\"{a}hler structure on $\mathcal{K}_{\Sigma}M$;

(2) For a normal disk $D\;$in $Map\left(  \Sigma,M\right)  _{emb}$, $\hat{D}$
is a holomorphic disk in $\mathcal{K}_{\Sigma}M$ if and only if $D\tilde
{\times}\Sigma$ gives a tight instanton in $M$ as above;

(3) $\mathcal{K}_{\Sigma}C$ is a Lagrangian submanifold in $\mathcal{K}%
_{\Sigma}M$ if and only if $C$ is a brane in $M$.
\end{theorem}

When we consider the complex analog of the above theorem on a Calabi-Yau
$n$-fold $M$, we choose any smooth manifold $\Sigma$ of dimension $n-2$, we
might hope that $Map\left(  \Sigma,M\right)  _{emb}/Diff\left(  \Sigma\right)
\otimes\mathbb{C}$, if exists, is hyperk\"{a}hler. Since the complexification
of $Diff\left(  \Sigma\right)  $ does not exist, we should interpret the above
quotient as a symplectic quotient $Map\left(  \Sigma,M\right)  _{emb}%
//Diff\left(  \Sigma\right)  $, if exists. The problem arises because one
needs to fix a background volume form on $\Sigma$ to define a symplectic
structure on $Map\left(  \Sigma,M\right)  _{emb}$. We will explain in the last
section on how to resolve this issue and prove the following theorem.

\begin{theorem}
Suppose $M$ is a Calabi-Yau $n$-fold and $\Sigma$ is a closed manifold of
dimension $n-2$. Then the isotropic knot space $\mathcal{\hat{K}}_{\Sigma}M$
has a natural holomorphic symplectic structure.
\end{theorem}

Furthermore the Calabi-Yau geometry on $M$ can be interpreted as the
holomorphic symplectic geometry on $\mathcal{\hat{K}}_{\Sigma}M$. Namely we
prove the following theorem in section \ref{Sec CLag Isot Knot}.

\begin{theorem}
Suppose $M$ is a Calabi-Yau $n$-fold. We have

(1) $\mathcal{\hat{K}}_{\Sigma}C$ is a $J$-complex Lagrangian submanifold in
$\mathcal{\hat{K}}_{\Sigma}M$ if and only if $C$ is a complex hypersurface in
$M$;

(2) $\mathcal{\hat{K}}_{\Sigma}C$ is a $K$-complex Lagrangian submanifold in
$\mathcal{\hat{K}}_{\Sigma}M$ if and only if $C$ is a special Lagrangian
submanifold in $M$ with phase $-\pi/2$.
\end{theorem}

Even though complex hypersurfaces and special Lagrangian submanifolds look
very different inside a Calabi-Yau manifold, their isotropic knot spaces are
both complex Lagrangian submanifolds in $\mathcal{\hat{K}}_{\Sigma}M$. One
reason is that any knot inside a special Lagrangian submanifold is
automatically isotropic. To prove this theorem, we need to construct carefully
certain appropiate deformations of isotropic knots inside $C$ so that
$\mathcal{\hat{K}}_{\Sigma}C$ being a $J$-complex Lagrangian (resp.
$I$-complex Lagrangian) implies that the dimension of $C$ is at least $2n-2$
(resp. at most $n$).

\section{Instantons and Branes}

In this section we introduce and study instantons and branes on manifolds with
(complex) vector cross products. Traditionally instantons refer to gradient
flow lines of a Morse function $f\,$on a Riemannian manifold $\left(
M,g\right)  $, as studied by Witten in \cite{W}. Morse theory can be
generalized to any closed one form $\phi$, because $\phi=df$ locally. Suppose
that $\phi$ is nonvanishing, we can choose a Riemannian metric on $M $ such
that it has unit length at every point. Then the gradient flow lines for the
vector field $X,$ defined by $\phi=\iota_{X}g$, can be reinterpreted as one
dimensional submanifolds in $M$ calibrated by $\phi$.

Suppose that $M$ is a symplectic manifold with symplectic form $\omega$. By
transgression on $\omega$, we obtain a closed one form on the free loop space
of $M$. The instantons in the free loop space correspond to holomorphic curves
in $M$ and they are calibrated by $\omega$. We continue to call these
instantons and they play important roles in the closed String theory (see e.g.
\cite{Mirror Mfd I}). In open String theory, we consider holomorphic curves in
$M$ with boundaries lying on a Lagrangian submanifold in $M$, which we call it
a brane.

In general instantons are submanifolds (of the smallest dimension) which are
preserved by the VCP, and branes are the natural boundaries for the free
boundary value problem for instantons.

\subsection{Vector Cross Product}

In this subsection we review the VCP as introduced by Gray. A basic example of
\textit{VCP} is the vector product$\;$in $\mathbb{R}^{3}$. The vector
product$\;$of any two linearly independent vectors in a plane is a vector
which is orthogonal to the plane and has the length equal to the area of the
parallelogram spanned by those vectors. Actually, these two properties
characterize \textit{VCP}$\;$as in the following definition by Brown and Gray
\cite{B-G}.

\begin{definition}
On an n-dimensional Riemannian manifold $M$ with a metric $g$, an r-fold
Vector Cross Product $\left(  \mathit{VCP}\right)  $ is a smooth bundle map,
\[
\chi:\wedge^{r}TM\rightarrow TM
\]
satisfying
\[
\left\{
\begin{tabular}
[c]{l}%
$g\left(  \chi\left(  v_{1},...,v_{r}\right)  ,v_{i}\right)  =0\;,\,\;\left(
1\leq i\leq r\right)  $\\
\\
$g\left(  \chi\left(  v_{1},...,v_{r}\right)  ,\chi\left(  v_{1}%
,...,v_{r}\right)  \right)  =\left\|  v_{1}\wedge...\wedge v_{r}\right\|
^{2}$%
\end{tabular}
\right.
\]
where $\left\|  \cdot\right\|  $ is the induced metric$\;$on\ $\wedge^{r}TM$.
\end{definition}

We will also denote
\[
v_{1}\times...\times v_{r}=\chi\left(  v_{1},...,v_{r}\right)  .
\]
The first condition in the above definition is equivalent to the following
tensor $\phi$ being a$\;$skew symmetric tensor of degree $r+1$, i.e. $\phi$ is
a differential form,
\[
\phi\left(  v_{1},...,v_{r+1}\right)  =g\left(  v_{1}\times...\times
v_{r},v_{r+1}\right)
\]

Clearly an $r$-fold VCP $\chi$ on $M$ can be characterized by an appropriate
differential form $\phi\in\Omega^{r+1}\left(  M\right)  $, which we call a
\textit{VCP form. }

\begin{definition}
An \textit{VCP} $form$\ of degree r+1 is a\ differential form $\phi$ $\in$
$\Omega^{r+1}\left(  M\right)  $\ satisfying
\[
\left|  \iota_{e_{1}\wedge e_{2}..\wedge e_{r}}(\phi)\right|  =1
\]
for any orthonormal $e_{1},...,e_{r}$ $\in$ $TM_{x}$, all $x$ in $M$.
\end{definition}

It is not difficult to see that a Hermitian almost complex structure is
equivalent to a $1$-fold \textit{VCP} and the corresponding K\"{a}hler form is
the corresponding \textit{VCP form}. In fact the complete list of VCPs is
surprisingly short. The classification of the linear \textit{VCP}'s on a
vector space $V$ \ with positive definite inner product $g$, by Brown and Gray
\cite{B-G}, can be summarized in the following.

(i) $r=1$ : If $\chi:V\rightarrow V$ is a $1$-fold VCP, then $\left|
\chi\left(  v\right)  \right|  =\left|  v\right|  $ implies that $\chi$ is an
orthogonal transformation. Polarizing $\left\langle \chi\left(  v\right)
,v\right\rangle =0$, we obtain $\left\langle \chi\left(  u\right)
,v\right\rangle +\left\langle u,\chi\left(  v\right)  \right\rangle =0$, that
is $\chi^{\ast}=-\chi$. Together we have $\chi^{2}=-id$. Namely a $1$-fold VCP
is equivalent to a Hermitian complex structure on $V$. The symmetry group of
$\chi$ is isomorphic to $U\left(  m\right)  =O\left(  2m\right)  \cap
GL\left(  m,\mathbb{C}\right)  $. On the other hand, the isomorphism $U\left(
m\right)  =O\left(  2m\right)  \cap Sp\left(  2m,\mathbb{R}\right)  $ reflects
the fact that a $1$-fold VCP is determined by its corresponding VCP form, or
its K\"{a}hler form $\phi$. The standard example is $V=\mathbb{C}^{m}$ with
\[
\phi=dx^{1}\wedge dy^{1}+\cdots+dx^{m}\wedge dy^{m}\text{.}%
\]

(ii) $r=n-1$ : an $\left(  n-1\right)  $-fold VCP on an $n$-dimensional inner
product space $V\;$is the Hodge star operator $\ast\;$given by $g$ on
$\Lambda^{n-1}V$ and the VCP form of degree $n$ is the induced volume form
$Vol_{V}$ on $V$. The automorphism group preserving the VCP form $Vol_{V}$ is
the group of linear transformations preserving $g$ and $Vol_{V}$, i.e.
$Aut\left(  V,Vol_{V}\right)  =O\left(  n\right)  \cap SL\left(
n,\mathbb{R}\right)  =SO\left(  n\right)  $. The standard example is
$V=\mathbb{R}^{n}$ and
\[
\phi=dx^{1}\wedge dx^{2}\wedge\cdots\wedge dx^{n}\text{.}%
\]

(iii) $r=2$ : a $2$-fold \textit{VCP on }a $7$-dimensional vector space
$\operatorname{Im}\mathbb{O}$ is a cross product defined as $a\times
b=\operatorname{Im}\left(  ab\right)  $ for any $a,b$ in $\operatorname{Im}%
\mathbb{O}$. For coordinates $\left(  x_{1,}...,x_{7}\right)  $ on
$\operatorname{Im}\mathbb{O}$, the corresponding \textit{VCP} form $\Omega$ of
degree $3$ can be written as
\[
\Omega=dx^{123}-dx^{167}+dx^{145}+dx^{257}+dx^{246}-dx^{356}+dx^{347}%
\]
where $dx^{ijk}=dx^{i}\wedge dx^{j}\wedge dx^{k}$. Bryant \cite{Br2} showed
that the group of real-linear transformations of $\operatorname{Im}\mathbb{O}$
preserving the \textit{VCP} form $\Omega$ actually preserves $g$ and VCP and
more it is exactly $G_{2}$, the automorphism group of the octonion
$\mathbb{O}$, i.e. $Aut\left(  \operatorname{Im}\mathbb{O,}\Omega\right)
=G_{2}\subset SO\left(  \operatorname{Im}\mathbb{O}\right)  =SO\left(
7\right)  $.

(iv) $r=$ $3$ : a $3$-fold \textit{VCP }on an $8$-dimensional vector space
$\mathbb{O}$ is a cross product defined as $a\times b\times c=$ $\frac{1}%
{2}\left(  a\left(  \bar{b}c\right)  -c\left(  \bar{b}a\right)  \right)  $ for
any $a,b$ and $c$ in $\mathbb{O}$. For coordinates $\left(  x_{1,}%
...,x_{8}\right)  $ on $\mathbb{O}$, the corresponding \textit{VCP} form
$\Theta$ of degree $4$ can be written as
\begin{align*}
\Theta &  =-dx^{1234}-dx^{5678}-\left(  dx^{21}+dx^{34}\right)  \left(
dx^{65}+dx^{78}\right) \\
&  -\left(  dx^{31}+dx^{42}\right)  \left(  dx^{75}+dx^{86}\right)  -\left(
dx^{41}+dx^{23}\right)  \left(  dx^{85}+dx^{67}\right)
\end{align*}
Bryant \cite{Br2} also showed that the group of real-linear transformations of
$\mathbb{O}$ preserving the \textit{VCP} form $\Theta$ on $\mathbb{O}$
preserves $g$ and VCP, and it is $Spin\left(  7\right)  $. i.e. $Aut\left(
\mathbb{O,}\Theta\right)  =Spin\left(  7\right)  \subset SO\left(  8\right)  $.%

\[
\underset{\text{Table 1: Classification of r-fold }\mathit{VCP}\;\text{on
V}\simeq\mathbb{R}^{n}}{%
\begin{tabular}
[c]{|c||c|c|c|c|}\hline
r$%
\begin{array}
[c]{l}%
\,\\
\,
\end{array}
$ & 1 & n-1 & 2 & 3\\\hline
$n=\dim V$ & 2m & n & 7$%
\begin{array}
[c]{l}%
\,\\
\,
\end{array}
$ & 8\\\hline
$V
\begin{array}
[c]{l}%
\,\\
\,
\end{array}
$ & $\mathbb{C}^{m}$ & $\mathbb{R}^{n}$ & Im$\mathbb{O}$ & $\mathbb{O}%
$\\\hline
\end{tabular}
}%
\]

From the above classification of linear VCPs, the existence of a VCP on a
Riemannian manifold $M$ is equivalent to the reduction of the structure group
of the frame bundle from $O\left(  n\right)  $ to $U\left(  m\right)  $,
$SO\left(  n\right)  ,$ $G_{2}$ and $Spin\left(  7\right)  $, for $r=1$,
$n-1$, $2 $ and $3$ respectively. Using the obstruction theory in topology,
the necessary and sufficient condition for the existence of such a reduction
of structure group in the $r=n-1$, $2$ and $3$ cases are $w_{1}=0$,
$w_{1}=w_{2}=0$ and $w_{1}=w_{2}=p_{1}^{2}-4p_{2}+8\chi=0$ respectively. Here
$w_{i}$, $p_{i}$ and $\chi$ are the Stiefel-Whitney class, the Pontrjagin
class, and the Euler class of $M$, respectively (see e.g. \cite{L-M}).

For a \textit{VCP} on a Riemannian manifold $M$, that is a linear VCP on each
tangent space of $M$, we would require certain integrability condition for its
coherence. For example, a $1$-fold VCP is a Hermitian almost complex structure
on $M$. This defines a symplectic structure or K\"{a}hler structure on $M$ if
the corresponding VCP form on $M$ is closed or parallel respectively. In
general we have the following definition.

\begin{definition}
Suppose that $M$ is a Riemannian manifold with a VCP $\chi$ and $\phi$ is its
corresponding VCP form. We call $\chi$ a closed (resp. parallel) VCP$\;$if
$d\phi=0$ (resp. $\nabla\phi=0$, where $\nabla$ is the Levi-Civita connection
on $M$).
\end{definition}

\bigskip

The classification of manifolds with closed/parallel VCP is presented in the
following table.
\[%
\begin{tabular}
[c]{|l|l|l|}\hline
$r
\begin{array}
[c]{l}%
\,\\
\,
\end{array}
$ & Closed VCP & Parallel VCP\\\hline\hline
$1
\begin{array}
[c]{l}%
\,\\
\,
\end{array}
$ & Almost K\"{a}hler manifolds & K\"{a}hler manifolds\\\hline
$n-1
\begin{array}
[c]{l}%
\,\\
\,
\end{array}
$ & Oriented manifolds & Oriented manifolds\\\hline
$2
\begin{array}
[c]{l}%
\,\\
\,
\end{array}
$ & Almost $G_{2}\text{-manifolds}$ & $G_{2}\text{-manifolds}$\\\hline
$3
\begin{array}
[c]{l}%
\,\\
\,
\end{array}
$ & $Spin\left(  7\right)  \text{-manifolds}$ & $Spin\left(  7\right)
\text{-manifolds}$\\\hline
\end{tabular}
\]

Remark: A \textit{VCP} form $\phi$ of degree $r+1$ on $M$ induces a
\textit{VCP} form of degree $r$ on any oriented hypersurface $H$ in $M$,
namely the restriction of $\iota_{\nu}\phi$ to $H$ where $\nu$ is the unit
normal vector field on $H$. However it usually does not preserve the
closedness of the VCP form. For example, Calabi \cite{C} and Gray \cite{G}
showed that such a two form on any hypersurface in $\mathbb{R}^{7}%
=\operatorname{Im}\mathbb{O}$ is never closed unless it is an affine
hyperplane. This result is generalized by Bryant \cite{Br1} to any codimension
two submanifold in $\mathbb{O}$, where he showed that closedness of such a two
form is equivalent to the submanifold being a complex hypersurfaces for some
suitable complex structure on $\mathbb{O}$.

\bigskip

Remark: We consider any nontrivial smooth map $f:M_{1}\rightarrow M_{2}$
between two Riemannian manifolds with $r$-fold VCPs. If $f$ preserves their
VCPs, i.e.
\[
f_{\ast}\left(  v_{1}\times...\times v_{r}\right)  =\left(  f_{\ast}%
v_{1}\right)  \times...\times\left(  f_{\ast}v_{r}\right)  ,
\]
for any tangent vectors $v_{i}$'s to $M_{1}$, then Gray \cite{G2} showed that
$f$ is actually an isometric immersion unless $r=1$ and $f$ would be a
holomorphic map. In particular any diffeomorphisms of $M$ preserving an
$r$-fold VCP with $r\geq2$ is automatically an isometry of $M$, i.e.
$Diff\left(  M,\chi\right)  \subset Isom\left(  M,g\right)  $.

In the next section we study submanifolds in $M$ which are preserved by the VCP.

\subsection{Instantons for Vector Cross Products}

In a symplectic manifold with a compatible almost complex structure $J$,
namely an almost K\"{a}hler manifold, an instanton is a two dimensional
submanifold which is preserved by $J$, and a brane is a Lagrangian
submanifold, i.e. a middle dimensional submanifold such that the restriction
of the symplectic form vanishes. They play essential roles in symplectic
geometry. These notions have natural analogs for manifolds with VCP.

\begin{definition}
Let $M$ be a Riemannian manifold with a closed r-$fold\;$\textit{VCP} $\chi$.
An $r+1$ dimensional submanifold $A$ is called an \textbf{instanton} if it is
preserved by $\chi$.
\end{definition}

From the definition, two instantons in $M$ can never intersect along a
codimension one subspace.

Suppose that $A$ is an instanton in $M$. At any given point $x$ in $A$ and any
tangent vectors $u_{1},...,u_{r-1}$ in $T_{x}A$ and any normal vector $\nu$ in
$N_{A/M,x}$, we have
\[
u_{1}\times...\times u_{r-1}\times\nu\in N_{A/M,x}.
\]
This is because for any $u$ in $T_{x}A$,
\begin{align*}
g\left(  u_{1}\times...\times u_{r-1}\times\nu,u\right)   &  =\phi\left(
u_{1},...,u_{r-1},\nu,u\right) \\
&  =-\phi\left(  u_{1},...,u_{r-1},u,\nu\right) \\
&  =-g\left(  u_{1}\times...\times u_{r-1}\times u,\nu\right) \\
&  =0.
\end{align*}
The last equality follows from the fact that $A$ is preserved by $\chi$ and
therefore $u_{1}\times...\times u_{r-1}\times u$ lies in $T_{x}A$. This
implies that there is bundle map,
\[
\Lambda^{r-1}T_{A}\otimes N_{A/M}\rightarrow N_{A/M}\text{.}%
\]
For example, when $r=1$, $A$ is a holomorphic curve in an almost K\"{a}hler
manifold $M$. Then the above bundle map is simply the complex structure of the
normal bundle of $A$ in $M$. In general the normal bundle to an instanton is
always a twisted spinor bundle over $A$ and the above bundle map is given by
the Clifford multiplication. This description should be useful in
understanding the deformation theory and moduli space of instantons.

\bigskip

An instanton in $M$ is always a minimal submanifold with absolute minimal
volume. This will follow from an equivalent characterization of instantons,
namely they are those submanifolds in $M$ calibrated by $\phi$. The theory of
calibration is developed by Harvey and Lawson in \cite{H-L} to produce
absolute minimal submanifolds. Recall that a closed differential form $\psi$
of degree $k$ on a Riemannian manifold $M$ is called a \emph{calibrating form
}if it satisfies
\[
\psi\left(  x\right)  |_{V}\leq Vol_{V},
\]
for every oriented $k$-plane $V\;$in $T_{x}M$, at each point $x$ in $M$. Here
$Vol_{V}$ is the volume form on $V$ for the induced metric. A \emph{calibrated
submanifold} $A\;$is a submanifold where $\psi|_{A}$ is equal to the induced
volume form on $A$. An important observation is the following: any other
submanifold $B$, homologous to $A$, satisfies
\[
Vol\left(  B\right)  \geq Vol\left(  A\right)  =\int_{A}\left[  \psi\right]
\text{,}%
\]
and the equality sign holds if and only if $B$ is also a calibrated
submanifold in $M$. As we will discover that all the examples of calibrated
submanifolds studied in \cite{H-L} are either instantons or branes in
manifolds with (complex) VCP.

The next lemma shows that instantons are calibrated.

\begin{lemma}
Let $M$ be a Riemannian manifold with a closed $r$-fold$\;$\textit{VCP} $\chi$
and we denote the corresponding \textit{VCP}$\;$form as $\phi$ . Then we have
(i) $\phi\;$is a calibrating form, (ii) an $\left(  r+1\right)  $-dimensional
submanifold $A$ in $M$ is calibrated by $\phi$ if and only if it is an instanton.
\end{lemma}

\begin{proof}
The closed form $\phi$ is a calibrating form because for any $x$ in $M$ and
for any oriented orthonormal tangent vectors $e_{1},e_{2},...,e_{r+1}$ at $x$,
it satisfies,%
\begin{align*}
\phi(e_{1},e_{2},...,e_{r+1})  &  =\langle\chi\left(  e_{1},e_{2}%
,...,e_{r}\right)  ,e_{r+1}\rangle\\
&  \leq\left|  \chi\left(  e_{1},e_{2},...,e_{r}\right)  \right|  \cdot\left|
e_{r+1}\right| \\
&  =\left\|  e_{1}\wedge e_{2}\wedge...\wedge e_{r}\right\|  \cdot\left|
e_{r+1}\right| \\
&  =1.
\end{align*}

The equality signs hold if and only if $\chi\left(  e_{1},e_{2},...,e_{r}%
\right)  =e_{r+1}$. Namely the linear span of $e_{i}$'s is preserved by $\chi$.
\end{proof}

\bigskip

As a corollary of this lemma and basic properties of calibration, An $\left(
r+1\right)  $-dimensional submanifold $A$ in $M$ is an instanton if and only
if its total volume with respect to the induced metric satisfies the following
equality,
\[
Vol\left(  A\right)  =\int_{A}\left[  \phi\right]  \text{.}%
\]

Remark: As a matter of fact, each term in the expansion of $\exp\left(
\phi\right)  $ is a calibrating form and the corresponding calibrated
submanifolds are preserved by the VCP $\chi$.

\textbf{Examples of instantons}: (i) Instantons is an oriented n-dimensional
manifold (i.e. VCP\ form of degree n) are open subsets in $M$. (ii) Instantons
in a K\"{a}hler manifold (i.e. parallel VCP\ form of degree 2) are holomorphic
curves. (iii) Instantons in a $G_{2}$-manifold (i.e. parallel VCP\ form of
degree 3) are called associative submanifolds. (iv) Instantons in a
$Spin\left(  7\right)  $-manifold (i.e. parallel VCP\ form of degree 4, called
the Cayley form) are called Cayley submanifolds (see \cite{H-L}).

Suppose that $\left(  M,g,\chi\right)  $ has a torus symmetry group, say
$M=X\times T^{k}$. Then an $r$-fold VCP on $M$ induces an $\left(  r-k\right)
$-fold VCP on $X$. Moreover a submanifold $B$ in $X$ is an instanton if and
only if $B\times T^{k}$ is an instanton in $M$. For instanton if $\Sigma$ is a
holomorphic curve is a Calabi-Yau threefold $X$, then $\Sigma\times S^{1}$ is
an associative submanifold in the $G_{2}$-manifold $X\times S^{1}$ and vice versa.

\begin{proposition}
Suppose that $\sigma$ is an isometric involution of $M$ preserving its
$r$-fold VCP $\chi$. If the fix point set $M^{\sigma}$ has dimension $r+1$
then it is an instanton in $M$.
\end{proposition}

Proof: Given any tangent vectors $v_{1},...,v_{r}$ to $M^{\sigma}$ at any
point $x$, we write
\[
v_{1}\times\cdots\times v_{r}=t+n,
\]
where $t$ (resp. $n$) is tangent (resp. normal) to $M^{\sigma}$. Since
$M^{\sigma}$ is the fix point set of an involution, we have $\sigma_{\ast
}\left(  t+n\right)  =t-n$. On the other hand,
\begin{align*}
\sigma_{\ast}\left(  v_{1}\times\cdots\times v_{r}\right)   &  =\sigma_{\ast
}\left(  v_{1}\right)  \times\cdots\times\sigma_{\ast}\left(  v_{r}\right) \\
&  =v_{1}\times\cdots\times v_{r}%
\end{align*}
because $\sigma$ preserves $\chi$. This implies $n=0$, thus $M^{\sigma}$ is
preserved by $\chi$, namely an instanton in $M$. $\blacksquare$

\bigskip

Remark: Given any $r$-dimensional analytic submanifold $S$ in $M$, an $r$-fold
VCP $\chi$ on $M$ determines a unique normal direction on $S$. Using
Cartan-K\"{a}hler theory, we can always integrate out this direction and
obtain an instanton in $M$ containing $S$ (see e.g. \cite{H-L}).

\bigskip

\subsection{\label{1Sec Deform Inst}Deformations of instantons}

In order to describe instantons and their deformations effectively, we need to
further develop the linear algebra of an inner product space $V$ with an
$r$-fold VCP
\[
\chi:\wedge^{r}V\rightarrow V\text{,}%
\]
or their associated VCP form $\phi\in\Lambda^{r+1}V^{\ast}$. We define
\[
\tau:\wedge^{r+1}V\rightarrow\wedge^{2}V
\]
as the composition of the following homomorphisms:
\[
\wedge^{r+1}V\overset{\text{(i)}}{\rightarrow}V\otimes\wedge^{r}%
V\overset{\text{(ii)}}{\rightarrow}V\otimes V\overset{\text{(iii)}%
}{\rightarrow}\wedge^{2}V,
\]
where these maps are (i) the natural inclusion, (ii) $id\otimes\chi$, (iii)
the natural projection. Explicitly we have
\[
\tau\left(  v_{1},...,v_{r+1}\right)  =\frac{1}{\sqrt{r+1}}\sum_{k=1}%
^{r+1}\left(  -1\right)  ^{k-1}v_{k}\wedge\chi\left(  v_{1},...,\widehat
{v_{k}},...,v_{r+1}\right)  \text{.}%
\]
As a matter of fact, the image of $\tau$ lies inside a much small subspace in
$\wedge^{2}V$. We define $\mathbf{g}\subset\mathbf{so}\left(  V\right)
\cong\wedge^{2}V^{\ast}$ to be the space of infinitesimal isometries of $V$
preserving the VCP $\chi$. Namely $\zeta\in\mathbf{so}\left(  V\right)
\subset End\left(  V\right)  $ lies inside $\mathbf{g}$ if
\[
\zeta\left(  v_{1}\times v_{2}\times...\times v_{r}\right)  =\zeta\left(
v_{1}\right)  \times v_{2}\times...\times v_{r}+...+v_{1}\times v_{2}%
\times...\times\zeta\left(  v_{r}\right)  \text{,}%
\]
for any $v_{i}$'s in $V$. This is equivalent to
\[
\sum_{i=1}^{r+1}\phi\left(  v_{1},v_{2},...,\zeta\left(  v_{i}\right)
,...,v_{r+1}\right)  =0,
\]
for any $v_{i}$'s in $V$. The next lemma says that the image of $\tau$ is
orthogonal to $\mathbf{g\subset}\Lambda^{2}V^{\ast}$ with respect to the
natural pairing between $\Lambda^{2}V$ and $\Lambda^{2}V^{\ast}$.

\begin{lemma}
Given any $r$-fold VCP $\chi$ on $V$, we have
\[
\tau:\wedge^{r+1}V\rightarrow\mathbf{g}^{\bot}\subset\wedge^{2}V.
\]
\end{lemma}

Proof: Given any $\zeta\in\mathbf{g\subset so}\left(  V\right)  \subset
End\left(  V\right)  $, we denote the corresponding two form in $\wedge
^{2}V^{\ast}$ as $\bar{\zeta}$. For any $v_{i}$'s in $V$, we compute
\begin{align*}
&  \left\langle \tau\left(  v_{1},...,v_{r+1}\right)  ,\bar{\zeta
}\right\rangle \\
&  =\bar{\zeta}\left(  \frac{1}{\sqrt{r+1}}\sum_{i=1}^{r+1}\left(  -1\right)
^{i-1}v_{i}\wedge\left(  v_{1}\times...\times\widehat{v_{i}}\times...\times
v_{r+1}\right)  \right) \\
&  =\frac{1}{\sqrt{r+1}}\left\langle \sum_{i=1}^{r+1}\left(  -1\right)
^{i-1}\left(  v_{1}\times...\times\widehat{v_{i}}\times...\times
v_{r+1}\right)  ,\zeta\left(  v_{i}\right)  \right\rangle \\
&  =\frac{\left(  -1\right)  ^{r}}{\sqrt{r+1}}\sum_{i=1}^{r+1}\phi\left(
v_{1},...,\zeta\left(  v_{i}\right)  ,...,v_{r+1}\right) \\
&  =0.
\end{align*}
Hence the result. $\blacksquare$

\bigskip

\begin{proposition}%
\[
\left|  \phi\left(  v_{1},...,v_{r+1}\right)  \right|  ^{2}+\left|
\tau\left(  v_{1},...,v_{r+1}\right)  \right|  _{\wedge^{2}}^{2}=\left|
v_{1}\wedge...\wedge v_{r+1}\right|  _{\Lambda^{r+1}}^{2}%
\]
\end{proposition}

Proof: It suffices to assume that $v_{i}$'s are orthogonal to each others.
Note that when $l\neq m$, we have $v_{l}$ perpendicular to both $v_{m}$ and
$\chi\left(  v_{1},...,\widehat{v_{m}},...,v_{r+1}\right)  $ and therefore,
\[
\left\langle v_{l}\wedge\chi\left(  v_{1},...,\widehat{v_{l}},...,v_{r+1}%
\right)  ,v_{m}\wedge\chi\left(  v_{1},...,\widehat{v_{m}},...,v_{r+1}\right)
\right\rangle =0\text{.}%
\]
We compute%

\begin{align*}
&  \left|  \tau\left(  v_{1},...,v_{r+1}\right)  \right|  ^{2}\\
&  =\frac{1}{r+1}\left|  \sum_{k=1}^{r+1}\left(  -1\right)  ^{k-1}v_{k}%
\wedge\chi\left(  v_{1},...,\widehat{v_{k}},...,v_{r+1}\right)  \right|
^{2}\\
&  =\frac{1}{r+1}\sum_{k=1}^{r+1}\left|  v_{k}\wedge\chi\left(  v_{1}%
,...,\widehat{v_{k}},...,v_{r+1}\right)  \right|  ^{2}\\
&  =\frac{1}{r+1}\sum_{k=1}^{r+1}\left(  \left|  v_{1}\right|  ^{2}%
\cdots\left|  v_{r+1}\right|  ^{2}-\left\langle v_{k},\chi\left(
v_{1},...,\widehat{v_{k}},...,v_{r+1}\right)  \right\rangle ^{2}\right) \\
&  =\frac{1}{r+1}\sum_{k=1}^{r+1}\left(  \left|  v_{1}\right|  ^{2}%
\cdots\left|  v_{r+1}\right|  ^{2}-\left|  \phi\left(  v_{1},...,v_{r+1}%
\right)  \right|  ^{2}\right) \\
&  =\left|  v_{1}\wedge...\wedge v_{r+1}\right|  ^{2}-\left|  \phi\left(
v_{1},...,v_{r+1}\right)  \right|  ^{2}\text{.}%
\end{align*}
Hence the result. $\blacksquare$

\bigskip

Remark: This proposition and the following corollary were first obtained by
Harvey and Lawson in the $G_{2}$- and $Spin\left(  7\right)  $-manifolds cases
using special structures of the octonions \cite{H-L}. As an immediately
corollory of the above proposition, we have

\begin{corollary}
\label{1Cor Lin Alg Instanton}Suppose that $\chi$ is an $r$-fold VCP on $V$.
An $\left(  r+1\right)  $-dimensional linear subspace $P\subset V$ is
preserved by $\chi$, i.e an instanton, if and only if
\[
\tau\left(  v_{1},...,v_{r+1}\right)  =0
\]
for any base $v_{i}$'s of $P$.
\end{corollary}

We also denote the above condition as $\tau\left(  P\right)  =0$. Assume this
is the case, we denote the orthogonal decomposition of $V$ as
\[
V=P\oplus N\text{.}%
\]
Similarly we have the following decomposition,
\[
\wedge^{2}V\cong\wedge^{2}P\oplus\wedge^{2}N\oplus\left(  P\otimes N\right)
\text{.}%
\]

\begin{proposition}
\label{1Prop Deform inst}Suppose that $\chi$ is an $r$-fold VCP on $V$ with an
orthogonal decomposition
\[
V=P\oplus N
\]
where $P$ is an $\left(  r+1\right)  $-dimensional linear subspace of $V$
preserved by $\chi$.\ We have
\[
\tau\left(  p_{1},...,p_{r},n\right)  \in P\otimes N
\]
for any $p_{i}$'s in $P$ and $n\in N$.
\end{proposition}

Proof: Note that $P$ being an instanton in $V$ implies that
\[
\left\langle \chi\left(  p_{1},...,p_{r-1},n\right)  ,p_{r}\right\rangle
=\pm\left\langle \chi\left(  p_{1},...,p_{r-1},p_{r}\right)  ,n\right\rangle
=0\text{.}%
\]
That is $\chi\left(  p_{1},...,p_{r-1},n\right)  \in N$, or equivalently,
$p_{r}\wedge\chi\left(  p_{1},...,p_{r-1},n\right)  \in P\otimes N$. However
$\tau\left(  p_{1},...,p_{r},n\right)  $ is a linear combination of terms of
this form and therefore $\tau\left(  p_{1},...,p_{r},n\right)  \in P\otimes
N$. Hence the proposition. $\blacksquare$

\bigskip

We are going to use these linear algebra results to study the deformations of
instantons. Given any VCP $\chi$ on a Riemannian manifold $M$, the spaces of
infinitesimal automorphisms of $\left(  T_{M},\chi\right)  $ on various fibers
glue together to form a subbundle
\[
\mathbf{g}_{M}\subset\Lambda^{2}T_{M}^{\ast}.
\]
Similarly the vector cross product $\chi$ determines a tensor%
\[
\tau\in\Omega^{r+1}\left(  M,\mathbf{g}_{M}^{\bot}\right)  \text{.}%
\]
On the global level, the above corollary \ref{1Cor Lin Alg Instanton} is
equivalent to the following result.

\begin{theorem}
If $M$ is a Riemannian manifold with an $r$-fold VCP $\chi$, then a $\left(
r+1\right)  $-dimensional submanifold $A$ is an instanton if and only if
\[
\tau|_{A}=0\in\Omega^{r+1}\left(  A,\mathbf{g}_{M}^{\perp}\right)  \text{.}%
\]
\end{theorem}

This theorem is useful in studying the deformation of instantons. We first
recall that any nearby submanifold to $A$ is the image of the exponential map%
\[
\exp_{v}:A\rightarrow M,
\]
for some small normal vector field $v\in\Gamma\left(  A,N_{A/M}\right)  $.
Therefore, given any instanton $A$ in $M$, we can describe its nearby
instantons as the zeros of the following map,%
\begin{align*}
F  &  :\Gamma\left(  A,N_{A/M}\right)  \rightarrow\Gamma\left(  A,\mathbf{g}%
_{M}^{\perp}\right) \\
\;F\left(  v\right)   &  =\ast_{A}\exp_{v}^{\ast}\left(  \tau\right)  ,
\end{align*}
where $\ast_{A}$ is the Hodge star operator of $A$. Suppose that $A_{t}$ is a
family of submanifolds in $M$ with $A=A_{0}$ an instanton and we denote its
variation normal vector field as%
\[
v=\left.  \frac{dA_{t}}{dt}\right|  _{t=0}\in\Gamma\left(  A,N_{A/M}\right)
\text{.}%
\]
By the proposition \ref{1Prop Deform inst} and $A$ being an instanton, we have%
\[
\left.  \frac{d\tau|_{A_{t}}}{dt}\right|  _{t=0}\in\Gamma\left(  A,T_{A}%
^{\ast}\otimes N_{A/M}\cap\mathbf{g}_{M}^{\perp}\right)  \text{.}%
\]
In particular the derivative of $F$ is given by%
\begin{gather*}
F^{\prime}\left(  0\right)  :\Gamma\left(  A,N_{A/M}\right)  \rightarrow
\Gamma\left(  A,T_{A}^{\ast}\otimes N_{A/M}\cap\mathbf{g}_{M}^{\perp}\right)
\\
F^{\prime}\left(  0\right)  \left(  \left.  \frac{dA_{t}}{dt}\right|
_{t=0}\right)  =\left.  \frac{d\tau|_{A_{t}}}{dt}\right|  _{t=0}\text{.}%
\end{gather*}

By studying individual cases, we find out that $N_{A/M}$ is always a twisted
spinor bundle over $A$ and the linearization of $F$ coincide with a twisted
Dirac operator, i.e. $F^{\prime}\left(  0\right)  =\mathcal{D}$.

\subsection{Branes for Vector Cross Products}

In symplectic geometry, the natural free boundary condition for holomorphic
curves require their boundaries lying inside a Lagrangian submanifold. The
analog of a Lagrangian submanifold for VCP forms of higher degree is called a brane.

\begin{definition}
Suppose $M$ is an n-dimensional manifold with a closed \textit{VCP}$\;$form
$\phi$ of degree r+1. A submanifold $C$ is called a \textbf{brane} if
\[
\left\{
\begin{tabular}
[c]{l}%
\begin{tabular}
[c]{l}%
$\phi\mid_{C}=0$%
\end{tabular}
$\;$\\
\\%
\begin{tabular}
[c]{l}%
$dimC=\left(  n+r-1\right)  /2.$%
\end{tabular}
\end{tabular}
\right.
\]
\end{definition}

Remark on the dimension of a brane : \label{Dimension of branes}Branes have
the largest possible dimension among submanifolds $C$ satisfying $\phi|_{C}=0
$. To see this, it suffices to consider the linear case. Take any $r-1$
dimensional linear subspace $W$ in $C$, the interior product of $\phi$ by any
orthonormal basis of $W$ determines a symplectic form on the orthogonal
complement of $W$ in $M$, which we denote as $M/W$. Furthermore $C/W$ is an
isotropic subspace in $M/W$ and therefore $\dim C/W\leq\left(  \dim
M/W\right)  /2$. The equality sign holds exactly when $\dim C=\left(
n+r-1\right)  /2$.

As we recall in a symplectic manifold, a holomorphic disk intersects
perpendicularly a Lagrangian submanifold along the boundary, we have the
following lemma for intersection of an instanton and brane along the boundary
of the instanton in a manifold with a closed VCP form.

\begin{lemma}
\label{Bounary}Let $A$ be an instanton in an $n$-dimensional manifold $M$ with
closed VCP form $\phi$ of degree $r+1.$ Suppose the boundary of $A$ lies in a
brane $C$, then $A$ intersect $C$ perpendicularly along $\partial A$.
\end{lemma}

\begin{proof}
For $x\in\partial A\subset C$, consider $u\in T_{x}A$ perpendicular to
$\partial A$ and any $v\in T_{x}C$. Observe that there are $u_{1},...,u_{r}\in
T_{x}\left(  \partial A\right)  $ such that $u=u_{1}\times u_{2}%
\times...\times u_{r}$ since $\phi|_{A}$ is the volume form on $A$. Then,
\begin{align*}
g\left(  u,v\right)   &  =g\left(  u_{1}\times u_{2}\times...\times
u_{r},v\right) \\
&  =\phi\left(  u_{1},u_{2},...,u_{r},v\right)  =0
\end{align*}
because $u_{i}$'s and $v$ lie in $T_{x}C$ and $\phi|_{C}=0$. That is $u$ is
perpendicular to $C$.
\end{proof}

\bigskip

Note that we only need the assumption $\phi|_{C}=0$ on $C$ in the above lemma.
The condition $\phi|_{C}=0$ also implies that $\left[  \phi\right]  \in
H^{r+1}\left(  M,C\right)  $. Any such instanton $A$ minimizes volume within
the relative homology class $\left[  A\right]  \in H_{r+1}\left(  M,C\right)
$ and with volume equals to the pairing of $\left[  \phi\right]  $ and
$\left[  A\right]  $. Furthermore any submanifold $\left(  A^{\prime},\partial
A^{\prime}\right)  \subset\left(  M,C\right)  $ with $\left[  A^{\prime
}\right]  =\left[  A\right]  $ and $vol\left(  A^{\prime}\right)  =vol\left(
A\right)  $ is also an instanton. However, if $\dim C<\left(  n+r-1\right)
/2$, then finding instantons with boundaries lying on $C$ is an overdetermined
system of equations.

\bigskip

Since the definition of a branes depends only on the closed VCP form $\phi$
instead of $\chi$, the image of any brane under an $\phi$-preserving
diffeomorphism $f\in Diff\left(  M,\phi\right)  $ is again a brane.
Infinitesimally, $v=df_{t}/dt|_{t=0}\in Vect\left(  M,\phi\right)  $ satisfies
$L_{v}\phi=0$. This implies that $\iota_{v}\phi$ is a closed form because
$\phi$ is closed.

\begin{definition}
Suppose that $\phi$ is a closed VCP form on $M$. An $\phi$-preserving vector
field $v\in Vect\left(  M,\phi\right)  $ is called an $\phi$%
\textbf{-Hamiltonian vector field} if $\iota_{v}\phi$ is exact. That is
\[
\iota_{v}\phi=d\eta
\]
for some degree $r$ differential form $\eta$, which we call an $\phi
$-\textbf{Hamiltonian differential form}.
\end{definition}

We will discuss the $\phi$-Hamiltonian equivalence of branes in the next section.

\bigskip

\textbf{Examples of branes}: (i) Branes in an oriented n-dimensional manifold
(i.e. VCP\ form of degree n) are hypersurfaces. (ii) Branes in a K\"{a}hler
manifold (i.e. parallel VCP\ form of degree 2) are Lagrangian submanifolds.
(iii) Branes in a $G_{2}$-manifold (i.e. parallel VCP\ form of degree 3) can
be identified as those four dimensional submanifolds calibrated by $\ast\phi$
(see \cite{H-L}) and they are called coassociative submanifolds. (iv) The next
proposition shows that there is no brane in any $Spin\left(  7\right)  $-manifold.

\begin{proposition}
\label{Brane in Spin(7)}Brane does not exist in any $Spin\left(  7\right)
$-manifold. That is there is no 5-dimensional submanifold where the Cayley
form vanishes.
\end{proposition}

\begin{proof}
Suppose that $C$ is any submanifold in a $Spin\left(  7\right)  $-manifold $M$
where the Cayley form vanishes. This implies that $\chi\left(  e_{i}%
,e_{j},e_{k}\right)  $ (denoted as $\chi_{ijk}$) is perpendicular to $C$ for
any orthonormal tangent vectors $e_{l}$'s on $C$. Notice that these unit
vectors satisfy
\[
\chi_{ijp}\,\bot\,\chi_{ijq}.
\]
This is because
\[
\Vert\chi(e_{i,}e_{j,}e_{p}+e_{q})\Vert=\Vert e_{p}+e_{q}\Vert,
\]
which implies that
\[
g\left(  \chi(e_{i,}e_{j},e_{p}),\chi(e_{i,}e_{j,}e_{q})\right)  =g\left(
e_{p},e_{q}\right)  =0.
\]

If $\dim C=5$, i.e. $C$ is a brane in $M$, then its normal bundle has rank
three. However, by the above property, $\chi_{123},\chi_{124},\chi_{134}\;$and
$\chi_{234}$ are four orthonormal vectors normal to $C$, which is a contradiction.
\end{proof}

\bigskip

\
\[
\underset{\text{Table\ 2:\ Classification of instantons and branes}}{%
\begin{tabular}
[c]{|c|c|c|c|}\hline
$%
\begin{array}
[c]{l}%
\text{Manifolds }M\,\\
\multicolumn{1}{c}{\,\text{(}\dim M\text{)}}%
\end{array}
$ & $%
\begin{array}
[c]{l}%
\text{\textit{VCP}\ form}\;\phi\,\\
\multicolumn{1}{c}{\,\text{(degree of }\phi\text{)}}%
\end{array}
$ & $%
\begin{array}
[c]{l}%
\text{Instanton\ }A\,\\
\multicolumn{1}{c}{\,\text{(}\dim A\text{)}}%
\end{array}
$ & $%
\begin{array}
[c]{l}%
\text{Brane}\;C\,\\
\multicolumn{1}{c}{\,\text{(}\dim C\text{)}}%
\end{array}
$\\\hline\hline
$%
\begin{array}
[c]{l}%
\text{Oriented mfd.}\,\\
\multicolumn{1}{c}{\,\text{(}n\text{)}}%
\end{array}
$ & $%
\begin{array}
[c]{l}%
\text{Volume form}\\
\multicolumn{1}{c}{\,\text{(}n\text{)}}%
\end{array}
$ &
\begin{tabular}
[c]{l}%
Open Subset\\
\multicolumn{1}{c}{($n$)}%
\end{tabular}
&
\begin{tabular}
[c]{l}%
Hypersurface\\
\multicolumn{1}{c}{($n-1$)}%
\end{tabular}
\\\hline
$%
\begin{array}
[c]{l}%
\text{K\"{a}hler mfd.}\,\\
\multicolumn{1}{c}{\,\text{(}2m\text{)}}%
\end{array}
$ & $%
\begin{array}
[c]{l}%
\text{K\"{a}hler form}\\
\multicolumn{1}{c}{\,\text{(}2\text{)}}%
\end{array}
$ &
\begin{tabular}
[c]{l}%
Holomorphic\\
Curve\\
\multicolumn{1}{c}{($2$)}%
\end{tabular}
&
\begin{tabular}
[c]{l}%
Lagrangian\\
Submanifold\\
\multicolumn{1}{c}{($m$)}%
\end{tabular}
\\\hline
$%
\begin{array}
[c]{l}%
\text{G}_{2}\text{-manifold}\,\\
\multicolumn{1}{c}{\,\text{(}7\text{)}}%
\end{array}
$ & $%
\begin{array}
[c]{l}%
\mathit{G}_{2}\text{-form}\\
\multicolumn{1}{c}{\,\text{(}3\text{)}}%
\end{array}
$ &
\begin{tabular}
[c]{l}%
Associative\\
Submanifold\\
\multicolumn{1}{c}{($3$)}%
\end{tabular}
&
\begin{tabular}
[c]{l}%
Coassociative\\
Submanifold\\
\multicolumn{1}{c}{($4$)}%
\end{tabular}
\\\hline
$%
\begin{array}
[c]{l}%
Spin\left(  7\right)  \text{-mfd.}\,\\
\multicolumn{1}{c}{\,\text{(}8\text{)}}%
\end{array}
$ & $%
\begin{array}
[c]{l}%
\text{Cayley form}\\
\multicolumn{1}{c}{\,\text{(}4\text{)}}%
\end{array}
$ &
\begin{tabular}
[c]{l}%
Cayley submfd.\\
\multicolumn{1}{c}{($4$)}%
\end{tabular}
& N/A\\\hline
\end{tabular}
}%
\]

Remark on \textit{0-fold }VCP: Even though we usually assume $r$ is positive
and exclude 0-fold \textit{VCP} in the classification, such a \textit{VCP }or
its corresponding VCP form is simply given by a closed one form $\phi$ with
unit pointwise length. When $\phi$ has integral period, we can integrate it to
obtain a function,
\[
f:M\rightarrow S^{1}\text{ and }\phi=f^{\ast}d\theta\text{.}%
\]
Instantons are gradient flow lines for the Morse function $f$ on $M$. Branes
are middle dimensional submanifolds in fibers of $f$.

\subsection{\label{1Sec Deform Brane}Deformation Theory of Branes}

The intersections theory of branes plays an important role in describing the
geometry of vector cross product, this is analoguous to the role of the
Floer's Lagrangian intersection theory in symplectic geometry. Deformation
theory of branes are essential in understanding both the intersections theory
of branes and the moduli space of branes.

First we need to identify the normal bundle to any brane. Note that $\phi
|_{C}=0$ implies
\[
\chi:\Lambda^{r}T_{C}\rightarrow N_{C/M}.
\]
When $C$ has the maximum possible dimension, i.e. a brane, this is a
surjective homomorphism onto $N_{C/M}$. It is because, otherwise, there exists
$\nu\in N_{C/M}$ perpendicular to the image of $\chi\left(  \Lambda^{r}%
T_{C}\right)  $, thus $\phi$ will vanish on the linear span of $T_{C}$ and
$\nu$, i.e. a bigger space containing $C$, a contradiction.

By taking the dual on $\chi:\Lambda^{r}T_{C}\rightarrow N_{C/M}$, we obtain an
injective map,
\[
t:N_{C/M}^{\ast}\rightarrow\Lambda^{r}T_{C}^{\ast}%
\]
defined by
\[
t\left(  \alpha\right)  \left(  u_{1},...,u_{r}\right)  :=\alpha\left(
\chi\left(  u_{1},...,u_{r}\right)  \right)
\]
for $\alpha\in N_{C/M}^{\ast}$ and $u_{1},...,u_{r}$ $\in$ $T_{C}$. Observe
\[
\alpha\left(  \chi\left(  u_{1},...,u_{r}\right)  \right)  =g\left(
\chi\left(  u_{1},...,u_{r}\right)  ,\breve{\alpha}\right)  =\phi\left(
u_{1},...,u_{r},\breve{\alpha}\right)
\]
where $\breve{\alpha}$ $\in N_{C/M}$ such that its dual is $\alpha$. Then, for
any $\alpha\in N_{C/M}^{\ast}$ with $\left|  \alpha\right|  =1,$ $t\left(
\alpha\right)  $ is a VCP form on $T_{C}$ of degree $r$. The reason is for any
orthonormal vector $e_{1},...,e_{r-1}\in$ $T_{C,x}$ and $x\in$ $C$,
\[
\left|  \iota_{e_{1}\wedge...\wedge e_{r-1}}\left(  t\left(  \alpha\right)
\right)  \right|  =\left|  -\iota_{_{e_{1}\wedge...\wedge e_{r-1}\wedge
\breve{\alpha}}}\left(  \phi\right)  \right|  =1
\]
because $\phi$ is a VCP form and $\breve{\alpha}$ is a unit normal vector.
Therefore we proved the following proposition.

\begin{proposition}
Suppose that $C$ is a brane in a manifold $M$ with a VCP\ form of degree
$r+1$, then the image of the map
\[
t:N_{C/M}\rightarrow\Lambda^{r}T_{C}^{\ast},
\]
is the subbundle spanned by VCP form of degree $r$ on $T_{C,x}$ for all $x\in
C$.
\end{proposition}

We denote $t\left(  N_{C/M}\right)  $ as $\Lambda_{VCP}^{r}T_{C}^{\ast}$.
Using in the classification of branes below, $\Lambda_{VCP}^{r}T_{C}^{\ast}$
equals to (i) $T_{C}^{\ast}$ when $r=1$, (ii) $\Lambda_{+}^{2}T_{C}^{\ast}$
when $r=2$ (iii) $\Lambda^{n-1}T_{C}^{\ast}$ when $r=n-1$. Note that brane
does not exist when $r=3$ (Proposition \ref{Brane in Spin(7)}).

Using the exponential map, small deformations of $C$ correspond to sections of
$N_{C/M}$ and the branes are the zeros of the following map,
\begin{gather*}
F:\Gamma\left(  N_{C/M}\right)  \rightarrow\Omega^{r+1}\left(  C\right) \\
F\left(  v\right)  =\left(  \exp_{v}\right)  ^{\ast}\phi\text{,}%
\end{gather*}
defined on a small neighborhood of the origin in $\Gamma\left(  N_{C/M}%
\right)  $. We are going to study the deformation theory of branes following
the approach by McLean in \cite{MC}. Under the identification $t:\Gamma\left(
N_{C/M}\right)  \cong\Omega_{VCP}^{r}\left(  C\right)  $, the differential of
$F$ at $0$ is given by the exterior derivative because
\[
dF\left(  0\right)  \left(  v\right)  =L_{v}\left(  \phi\right)
|_{C}=d\left(  \iota_{v}\phi\right)  |_{C}=d\left(  t\left(  v\right)
\right)  \text{.}%
\]

Recall $F\left(  0\right)  =0$, we obtain $\left[  F\left(  v\right)  \right]
=\left[  F\left(  0\right)  \right]  =0\in H^{r+1`}\left(  C\right)  $ because
$C$ and $\exp_{v}\left(  C\right)  $ are homologous in $M$. Therefore we have
\begin{gather*}
F:\Omega_{VCP}^{r}\left(  C\right)  \rightarrow d\Omega^{r}\left(  C\right) \\
F\left(  0\right)  =0,\\
dF\left(  0\right)  =d\text{.}%
\end{gather*}
If we know that $d\Omega_{VCP}^{r}\left(  C\right)  =d\Omega^{r}\left(
C\right)  $, then using the implicit function theorem, we can show that
$F^{-1}\left(  0\right)  $ is smooth near $0$ and the tangent space is given
by the kernel of $dF\left(  0\right)  $. The condition $d\Omega_{VCP}%
^{r}\left(  C\right)  =d\Omega^{r}\left(  C\right)  $ can be verified in each
individual case, however the authors do not know of any general proof of this.
In any case we have proved the following result.

\begin{proposition}
Suppose that $\phi$ is a VCP form of degree $r+1$ on $M$. Then small
deformations of any brane $C$ are parametrized by closed form in $\Omega
_{VCP}^{r}\left(  C\right)  $. In particular the space of branes in $M$ is smooth.
\end{proposition}

\bigskip

The space of branes is usually of infinite dimensional. But quotienting out
the equivalence relationship of $\phi$-Hamlitonian, the moduli space of branes
of finite dimensional.

\begin{definition}
Suppose that $C_{1}$ and $C_{0}$ are two branes in a manifold $M$ with a
closed VCP form $\phi$ of degree $r+1$. They are called $\phi$%
\textbf{-Hamlitonian equivalent} to each other if they are joined by a family
of branes $C_{t}$ such that their deformation vector fields $v_{t}%
=dC_{t}/dt\in\Gamma\left(  N_{C_{t}/M}\right)  $ satisfy
\[
\iota_{v_{t}}\phi=d\eta_{t},
\]
for some $\eta_{t}\in\Omega^{r-1}\left(  C\right)  $.
\end{definition}

Using the Hodge theory and the previous proposition, we have the tangent space
to the moduli space of branes at any point $C$ equals $H_{VCP}^{r}\left(
C\right)  $, the space of harmonic forms in $\Omega_{VCP}^{r}\left(  C\right)
$. In particular, the moduli space is smooth and of finite dimensional. Its
tangent space is given by (i) $H^{1}\left(  C,\mathbb{R}\right)  $ when $r=1$;
(ii) $H_{+}^{2}\left(  C,\mathbb{R}\right)  $ when $r=2$ and (iii)
$H^{n-1}\left(  C,\mathbb{R}\right)  \cong\mathbb{R}$ when $r=n-1$. In the
third case, i.e. $\phi$ is the volume form on $M$, two nearby hypersurfaces
$C$ and $C^{\prime}$ are $\phi$-Hamiltonian equivalent if there is a singular
chain $B$ satisfying $\partial B=C-C^{\prime}$ and $Vol\left(  B\right)  =0$.

\bigskip

Lagrangian intersection theory in symplectic geometry plays the central role
in the subject, and also plays very essential roles in mirror symmetry.
Naively speaking we need to \textit{count} the number of instantons bounding
two Lagrangian submanifolds. It is natural to generalize this to other VCPs
and count the number of instantons bounding two branes. This is a very
difficult problem except when $r$ equals zero. In this case, suppose that
$C_{1}$ and $C_{2}$ are two branes in $M$, i.e. $C_{i}\subset f^{-1}\left(
\theta_{i}\right)  $ for $i=1,2$ are middle dimensional submanifolds. Here we
continue the notations in the previous remark. Since $f$ is a Riemannian
submersion, $M$ is a Riemannian mapping cylinder, i.e.
\[
M=X\times\left[  0,1\right]  /\sim
\]
for some isometry $h$ on $X$, identifying $X\times\left\{  0\right\}  $ and
$X\times\left\{  1\right\}  $. Thus both $C_{i}$'s can be regarded as middle
dimensional submanifolds in $X$. Then instantons in $M$ bounding $C_{1}$ and
$C_{2}$ correspond to intersection points between $C_{1}$ and $h^{k}\left(
C_{2}\right)  $ in $X$ for any integer $k$. Therefore the generating function
for the number of instantons is given explicitly by the following topological
sum,
\[
\sum_{k=-\infty}^{\infty}\#\left(  C_{1}\cap h^{k}\left(  C_{2}\right)
\right)  t^{k}.
\]

\section{Complex Vector Cross Product}

\subsection{\label{CVCP}Classification of Complex Vector Cross Products}

In this section we study vector cross products on complex vector spaces, or
Hermitian complex manifolds. For the sake of convenience, the complex vector
cross product ($\mathbb{C}$-\textit{VCP})\ will be defined in terms of complex
vector cross forms on a Hermitian complex manifold. Recall that a Hermitian
complex manifold is a Riemannian manifold $\left(  M,g\right)  $ with a
Hermitian complex structure $J$, that is $g\left(  Ju,Jv\right)  =g\left(
u,v\right)  $ for any tangent vectors $u$ and $v$.

\begin{definition}
\label{Complex VCP}On\ a Hermitian complex manifold $\left(  M,g,J\right)  $
of complex dimension n, an r-fold \textbf{complex vector cross product}
(abbrev. $\mathbb{C}$-$\mathit{VCP}$) is a holomorphic form $\phi$ of degree
$r+1$ satisfying
\[
\left|  \iota_{e_{1}\wedge e_{2}..\wedge e_{r}}(\phi)\right|  =2^{\left(
r+1\right)  /2},
\]
for any orthonormal tangent vectors $e_{1},...,e_{r}$ $\in$ $T_{x}^{1,0}M$,
for any $x$ in $M$.

A $\mathbb{C}$-VCP is called closed (resp. parallel) if $\phi$ is closed
(resp. parallel with respect to the Levi-Civita connection) form.
\end{definition}

Notice that if the manifold $M$ is a closed K\"{a}hler manifold, then every
holomorphic form in $M$ is closed. For completeness we include the proof of
this well-known fact.

\begin{lemma}
Suppose $M$ is a closed K\"{a}hler manifold. Then every holomorphic form is a
closed differential form.
\end{lemma}

\begin{proof}
Assume that $\psi$ is any holomorphic form of degree $k$ in $M$, that is
$\psi\in\Omega^{k,0}\left(  M\right)  $ and $\bar{\partial}\psi=0$. We need to
show that $\partial\psi=0\in\Omega^{k+1,0}\left(  M\right)  $. By the Riemann
bilinear relation, the pairing
\[
\int\eta_{1}\wedge\bar{\eta}_{2}\wedge\omega^{n-k-1}%
\]
is definite on $\eta_{i}\in\Omega^{k+1,0}\left(  M\right)  $. That is,
\[
\int_{M}\left|  \partial\psi\right|  ^{2}\omega^{n}=C\int_{M}\partial
\psi\wedge\overline{\partial\psi}\wedge\omega^{n-k-1}.
\]
Using integration by part on closed manifolds and holomorphicity of $\psi$, we
have
\[
\int_{M}\left|  \partial\psi\right|  ^{2}\omega^{n}=-C\int_{M}\partial
\bar{\partial}\psi\wedge\bar{\psi}\wedge\omega^{n-k-1}=0.
\]
This implies that $\partial\psi=0$, that is $\psi$ is a closed form on $M$.
\end{proof}

\bigskip

We are going to see that there are exactly two classes of K\"{a}hler manifolds
with $\mathbb{C}$-VCP, namely Calabi-Yau manifolds and hyperk\"{a}hler
manifolds. Furthermore every $\mathbb{C}$-VCP is automatically parallel, in
particular closed, provided that the manifold itself is closed.

\bigskip

\textbf{Example}: \textit{Calabi-Yau manifold} (i.e. $(n-1)$-fold $\mathbb{C}%
$-VCP). A linear complex volume form $\phi$ on $\mathbb{C}^{n}$ is an element
in $\Lambda^{n,0}\left(  \mathbb{C}^{n}\right)  $ with $\phi\bar{\phi}$ equals
the Riemannian volume form on $\mathbb{C}^{n}\cong\mathbb{R}^{2n}$. This is
because of the equality $\left|  \det_{\mathbb{C}}\left(  A\right)  \right|
^{2}=\det\left(  A_{\mathbb{R}}\right)  $ between a complex matrix $A$ and its
realization $A_{\mathbb{R}}$. It is given as follow,
\[
\phi=dz^{1}\wedge dz^{2}\wedge...\wedge dz^{n}\text{,}%
\]
for a suitable choice of complex coordinate $z^{j}$'s on $V$. It is easy to
see that $\phi$ defines a constant $\left(  n-1\right)  $-fold $\mathbb{C}%
$-VCP. Similarly an $\left(  n-1\right)  $-fold $\mathbb{C}$-VCP structure on
a closed K\"{a}hler manifold $\left(  M,g\right)  $ is a holomorphic volume
form $\Omega\in\Omega^{n,0}\left(  M\right)  $,
\[
\Omega\bar{\Omega}=C_{n}\omega^{n},
\]
where the constant $C_{n}$ equals $i^{n}\left(  -1\right)  ^{n\left(
n-1\right)  /2}2^{-n}/n!.$ This implies that the Ricci curvature of $M$
vanishes. Thus, using Bochner arguments, we can show that every holomorphic
form on $M$ is parallel. In particular $\Omega$ is a parallel complex volume
form on $M$ and therefore the holonomy group of $M$ lies inside $SU\left(
n\right)  $, i.e. a Calabi-Yau manifold. A celebrated theorem of Yau \cite{Yau
Calabi Conj AlgGeom} says that any closed K\"{a}hler manifold with trivial
first Chern class $c_{1}\left(  M\right)  $ admits K\"{a}hler metric with
vanishing Ricci curvature. This implies that $M$ is a Calabi-Yau manifold if
the canonical line bundle is trivial holomorphically.

\bigskip

\textbf{Example}: \textit{Hyperk\"{a}hler manifold} (i.e. 1-fold $C$-VCP). A
hyperk\"{a}hler manifold is a Riemannian manifold $\left(  M,g\right)  $ of
dimension $n=4m$ with its holonomy group lies inside $SU\left(  m\right)
=GL\left(  m,\mathbb{H}\right)  \cap SO\left(  4m\right)  $. Namely it has
parallel Hermitian complex structures $I$, $J$ and $K$ satisfying the Hamilton
relation,
\[
I^{2}=J^{2}=K^{2}=IJK=-Id\text{.}%
\]
These complex structures defines three different K\"{a}hler structures
$\omega_{I}$, $\omega_{J}$ and $\omega_{K}$ on $\left(  M,g\right)  $
respectively. If we fix one of them, say $J$, then $\Omega=\omega_{I}%
-i\omega_{K}\in\Omega^{2,0}\left(  M\right)  $ is a parallel $J$-holomorphic
symplectic form on $M$. These two descriptions of a hyperk\"{a}hler manifold
are equivalent and it is simply the global version of the isomorphism
$Sp\left(  m\right)  =U\left(  2m\right)  \cap Sp\left(  m,\mathbb{C}\right)
$. This form $\Omega$ is a parallel $1$-fold $\mathbb{C}$-VCP form on $M$. The
reasoning is the same as the one in the real case. In the linear case, this is
given as follow,%

\[
\Omega=dz^{1}\wedge dz^{2}+......+dz^{2m-1}\wedge dz^{2m}\text{,}%
\]
for some suitable choice of coordinates on $\mathbb{C}^{2m}$. Conversely, if
$\Omega$ is a $1$-fold $\mathbb{C}$-VCP on a closed K\"{a}hler manifold
$\left(  M,g,J\right)  $, then it is a holomorphic symplectic form on $M$.
Since $Sp\left(  m\right)  \subset SU\left(  2m\right)  $, any hyperk\"{a}hler
manifold is a Calabi-Yau manifold. This implies that $\Omega$ is indeed
parallel as before. Therefore a hyperk\"{a}hler structure is equivalent to a
$1$-fold $\mathbb{C}$-VCP on any closed K\"{a}hler manifold.

\bigskip

We remark that, as in the real setting, a constant $r$-fold $\mathbb{C}$-VCP
on a complex vector space induces an $\left(  r-1\right)  $-fold $\mathbb{C}%
$-VCP on any of its complex hyperplane.

We are going to show that there is no other complex vector cross product
besides the holomorphic volume form and the holomorphic symplectic form as
discussed above. In particular, there is no complex analog of \textit{VCP} for
$G_{2}$-manifolds and $Spin\left(  7\right)  $-manifolds.

\begin{proposition}
On a complex vector space $V\;$of complex dimension $n$, there is an $r$-fold
$\mathbb{C}$-VCP if and only if either (i) $r$ $=1$ and $n=2m$ or (ii) $r=n-1$
and $n$ arbitrary. The corresponding $\mathbb{C}$-VCP form is a holomorphic
symplectic form and a holomorphic volume form respectively
\end{proposition}

\begin{proof}
From the above two examples, there is an (n-1)-fold$\;\mathbb{C}$-\textit{VCP}
on $V$, and more if n is an even number, a 1-fold $\mathbb{C}$-\textit{VCP}
exists on it. Now, we need to see there is no other type of $\mathbb{C}%
$-\textit{VCP} on a complex vector space. For that matter, we claim that for
$r\geq2$ if there is an $r$-fold$\;\mathbb{C}$-\textit{VCP} on a complex
vector space $V\;$of complex dimension n, $r$ must be n-1. At first, observe
that for $r\geq2$, an $r$-fold $\mathbb{C}$-\textit{VCP} on a vector space
induces an $\left(  r-1\right)  $-fold $\mathbb{C}$-\textit{VCP} on the
complex hyperplane. Therefore, an $r$-fold $\mathbb{C}$-\textit{VCP} on a
complex vector space of complex dimension $n$ is reduced $2$-fold $\mathbb{C}%
$-\textit{VCP} on a complex $(n-r+2)$-dimensional vector space. Now, to show
claim, it is enough to verify if there is a $2$-fold $\mathbb{C}$-\textit{VCP}
on a complex vector space $W$, then its complex dimension must be $3$.

As in the one of the example of $\mathbb{C}$-\textit{VCP}, when
dim$_{\mathbb{C}}$ $W$ $=3$, there is a $2$-fold $\mathbb{C}$-\textit{VCP}.
Now, we need to show there is no higher complex vector space with $2$-fold
$\mathbb{C}$-\textit{VCP}. Suppose dim$_{\mathbb{C}}$ $W$ $\geq4$ with 2-fold
$\mathbb{C}$-\textit{VCP}$\;\phi$, and by choosing any unit holomorphic vector
$z$ in $W$, consider a complex subspace $Z\;$spanned by $z$ and $\bar{z}$.
Then,$\;\iota_{z}\left(  \phi\right)  $ is a $1$-fold $\mathbb{C}%
$-\textit{VCP} on $Z^{\perp}$ and more $\dim_{\mathbb{C}}\left(  Z^{\perp
}\right)  $ is at least $4$ because $\dim_{\mathbb{C}}\left(  Z^{\perp
}\right)  \geq3$ and an even number so that $Z^{\perp}$ has a 1-fold
$\mathbb{C}$-\textit{VCP}.

Now, we may rewrite the 2-fold\ $\mathbb{C}$-\textit{VCP}$\;$on $W$ as
\[
\phi=z^{\ast}\wedge\iota_{z}\left(  \phi\right)  +\phi_{1},
\]
where $z^{\ast}$ is dual form of $z$ and $\phi_{1}$ is the sum of
terms\ without\ $z^{\ast}$. Since $\dim_{\mathbb{C}}\left(  Z^{\perp}\right)
$ is at least 4, 1-fold$\;\mathbb{C}$-\textit{VCP}$,\;\iota_{z}\left(
\phi\right)  $ on $Z^{\perp}$ has of the form $a_{1}^{\ast}\wedge b_{1}^{\ast
}+a_{2}^{\ast}\wedge b_{2}^{\ast}....$where $a_{1},a_{2},b_{1}$ and $b_{2}$
are orthonormal holomorphic vectors in $Z^{\perp}$. We consider the following,%
\begin{align*}
\iota_{\left(  b_{1}+b_{2}\right)  }\iota_{\left(  a_{1}+a_{2}\right)
}\left(  \phi\right)   &  =\iota_{\left(  b_{1}\right)  }\iota_{\left(
a_{1}\right)  }\left(  \phi\right)  +\iota_{\left(  b_{2}\right)  }%
\iota_{\left(  a_{1}\right)  }\left(  \phi\right)  +\iota_{\left(
b_{1}\right)  }\iota_{\left(  a_{2}\right)  }\left(  \phi\right)
+\iota_{\left(  b_{2}\right)  }\iota_{\left(  a_{2}\right)  }\left(
\phi\right) \\
&  =2z^{\ast}+\iota_{\left(  b_{2}\right)  }\iota_{\left(  a_{1}\right)
}\left(  \phi\right)  +\iota_{\left(  b_{1}\right)  }\iota_{\left(
a_{2}\right)  }\left(  \phi\right)  ,
\end{align*}
and
\[
\iota_{\left(  -\sqrt{-1}b_{1}+b_{2}\right)  }\iota_{\left(  -\sqrt{-1}%
a_{1}+a_{2}\right)  }\left(  \phi\right)  =-\sqrt{-1}\iota_{\left(
b_{2}\right)  }\iota_{\left(  a_{1}\right)  }\left(  \phi\right)  -\sqrt
{-1}\iota_{\left(  b_{1}\right)  }\iota_{\left(  a_{2}\right)  }\left(
\phi\right)  .
\]
Note that $a_{1}+a_{2}$, $b_{1}+b_{2}$, $-\sqrt{-1}a_{1}+a_{2}\;$and
$-\sqrt{-1}b_{1}+b_{2}$ are holomorphic vectors with the same length.

The interior product of any term from $\iota_{z}\left(  \phi\right)  $ and
$\phi_{1}$ is zero because it satisfies for example , $\left|  \iota
_{a_{1}\wedge b_{1}}\left(  \phi\right)  \right|  =1$. This implies that
$\phi_{1}$ does not have any term in $\iota_{z}\left(  \phi\right)  $. Hence
we have $\iota_{b_{1}}\iota_{a_{1}}\left(  \phi\right)  $ = $z^{\ast}$ and
$\iota_{b_{2}}\iota_{a_{2}}\left(  \phi\right)  =$ $z^{\ast}$.

From the choice of all orthonormal $z,a_{1},b_{1},a_{2}$ and $b_{2}$,
holomorphic vectors $a_{1}+a_{2}$ and $b_{1}+b_{2}$ are orthogonal to each
other, and this is true between holomorphic vectors $-\sqrt{-1}a_{1}+a_{2}%
\;$and $-\sqrt{-1}b_{1}+b_{2}$. So $\iota_{\left(  b_{1}+b_{2}\right)  }%
\iota_{\left(  a_{1}+a_{2}\right)  }\left(  \phi\right)  $ and $\iota_{\left(
-\sqrt{-1}b_{1}+b_{2}\right)  }\iota_{\left(  -\sqrt{-1}a_{1}+a_{2}\right)
}\left(  \phi\right)  $ are supposed to produce the same length by definition
of $\mathbb{C}$-\textit{VCP}, but it can be checked that is impossible because
$z^{\ast}$ and $\iota_{\left(  b_{2}\right)  }\iota_{\left(  a_{1}\right)
}\left(  \phi\right)  +\iota_{\left(  b_{1}\right)  }\iota_{\left(
a_{2}\right)  }\left(  \phi\right)  $ are perpendicular each other. From this
contradiction, a complex vector space $W$ is of complex dimension 3 so that it
has a 2-fold\ $\mathbb{C}$-\textit{VCP}.
\end{proof}

From this proposition and the examples of $\mathbb{C}$-\textit{VCP}, one can
conclude the following theorem.

\begin{theorem}
\label{Classification of CVCP}(\textbf{Classification of }$\mathbb{C}%
$\textbf{-}\textit{VCP}) Suppose $M$ is a closed K\"{a}hler manifold of
complex dimension n, with an r -fold $\mathbb{C}$-\textit{VCP}. Then either

(i) $r=n-1$ and $M$ is a Calabi-Yau manifold, or

(ii) $r=1$ and $M$ is a hyperk\"{a}hler manifold.
\end{theorem}

\subsection{Instantons for Complex Vector Cross Products}

In this section, we introduce and study instantons and branes on a K\"{a}hler
manifold $M$ with a closed $\mathbb{C}$-VCP $\phi\in\Omega^{r+1,0}\left(
M\right)  $. Recall that an instanton, in the real setting, is a $\left(
r+1\right)  $-dimensional submanifold $A$ preserved by $\chi$, or equivalently
$A$ is calibrated by the VCP form. In the complex setting, the real and
imaginary parts of the complex VCP form are always calibrating forms and we
called such calibrated submanifolds instantons.

\begin{lemma}
Suppose $\phi$ is a closed $\mathbb{C}$-VCP form of degree $r+1$ on a
K\"{a}hler manifold $M$, then (i) $\operatorname{Re}\left(  e^{i\theta}%
\phi\right)  $ is a calibrating form for any real number $\theta$, and (ii) a
$\left(  r+1\right)  $-dimensional submanifold $A$ in $M$ is calibrated by
$\operatorname{Re}\left(  e^{i\theta}\phi\right)  $ only if
\[
\operatorname{Im}\left(  e^{i\theta}\phi\right)  |_{A}=\text{ }\omega
|_{A}=0\text{.}%
\]
\end{lemma}

\begin{proof}
It suffices to check the linear case, namely $M=\mathbb{C}^{n}$ with the
standard complex structure $J.$ Consider any oriented orthonormal vectors
$a_{1}$,..., $a_{r+1}$ in $\mathbb{R}^{2n}=\mathbb{C}^{n}$ and denote
$\xi=a_{1}\wedge...\wedge$ $a_{r+1},$then
\[
\left\{  \operatorname{Re}\left(  e^{i\theta}\phi\right)  \left(  \xi\right)
\right\}  ^{2}+\left\{  \operatorname{Im}\left(  e^{i\theta}\phi\right)
\left(  \xi\right)  \right\}  ^{2}=\left|  \left(  e^{i\theta}\phi\right)
\left(  \xi\right)  \right|  ^{2}=\left|  \phi\left(  \xi\right)  \right|
^{2}.
\]
Since $\phi$ is of type $\left(  r+1,0\right)  $, we have
\[
\left|  \phi\left(  \xi\right)  \right|  ^{2}=2^{-\left(  r+1\right)  }\left|
\phi\left(  \xi\otimes\mathbb{C}\right)  \right|  ^{2}.
\]
where $\xi\otimes\mathbb{C}=\tilde{a}_{1}\wedge...\wedge$ $\tilde{a}_{r+1}$,
with $\tilde{a}_{i}\mathbb{=:}\left(  a_{i}-\sqrt{-1}Ja_{i}\right)  /\sqrt
{2}.$ This is because $dz_{i}\left(  \tilde{a}_{k}\right)  =\sqrt{2}%
dz_{i}\left(  a_{k}\right)  $.

If we denote the dual vector of any one form $\eta$ as $\eta^{\#}$, then%

\begin{align*}
2^{-\left(  r+1\right)  }\left|  \phi\left(  \xi\otimes\mathbb{C}\right)
\right|  ^{2}  &  =2^{-\left(  r+1\right)  }\left|  \phi\left(  \tilde{a}%
_{1},...,\tilde{a}_{r+1}\right)  \right|  ^{2}\\
&  =2^{-\left(  r+1\right)  }\left|  \left\langle \left(  \iota_{\tilde{a}%
_{1}\wedge...\wedge\tilde{a}_{r}}\phi\right)  ^{\#},\tilde{a}_{r+1}%
\right\rangle \right| \\
&  \leq2^{-\left(  r+1\right)  }\left|  \left(  \iota_{\tilde{a}_{1}%
\wedge...\wedge\tilde{a}_{r}}\phi\right)  ^{\#}\right|  \left|  \tilde
{a}_{r+1}\right| \\
&  =2^{-\left(  r+1\right)  }\left|  \iota_{\tilde{a}_{1}\wedge...\wedge
\tilde{a}_{r}}\phi\right|  \leq1,
\end{align*}
because $\phi$ is a $\mathbb{C}$-VCP and $\tilde{a}_{i}^{\prime}$s are
elements in $T^{1,0}M$ of unit length. Therefore
\[
\left|  \operatorname{Re}\left(  e^{i\theta}\phi\right)  \left(  \xi\right)
\right|  \leq1\text{,}%
\]
and when the equality sign holds, we have (i) $\left|  \operatorname{Im}%
\left(  e^{i\theta}\phi\right)  \left(  \xi\right)  \right|  =0$, (ii)
$\left(  \iota_{\tilde{a}_{1}\wedge...\wedge\tilde{a}_{r}}\phi\right)  ^{\#}$
parallel to $\tilde{a}_{r+1}$ and (iii) $\left|  \iota_{\tilde{a}_{1}%
\wedge...\wedge\tilde{a}_{r}}\phi\right|  =2^{r+1}$. Since $a_{i}$'s are
orthonormal, condition (iii) is equivalent to the $\tilde{a}_{i}$'s being
orthonormal vectors in $T^{1,0}M$. This happens exactly when the linear span
of $a_{i}$'s is isotropic with respect to $\omega$.
\end{proof}

We remark that when $\phi$ is the holomorphic volume form of a Calabi-Yau
manifold, then a middle dimensional submanifold $A$ in $M$ is calibrated by
$\operatorname{Re}\left(  e^{i\theta}\phi\right)  $ if and only if it
satisfies $\operatorname{Im}\left(  e^{i\theta}\phi\right)  |_{A}=$
$\omega|_{A}=0$, and it is called a special Lagrangian submanifold with phase
angle $\theta$.

\begin{definition}
On a closed K\"{a}hler manifold $M$ with an r-fold $\mathbb{C}$-\textit{VCP}
$\phi$, an r +1 dimensional submanifold $A$ is called an\textbf{\ instanton}
with phase $\theta$\ $\in\mathbb{R}\;$if it is calibrated by $Re\left(
e^{i\theta}\phi\right)  $ $,$ i.e.
\[
Re\left(  e^{i\theta}\phi\right)  \mid_{A}=vol_{A}\text{.}%
\]
Equivalently, $\operatorname{Im}\left(  e^{i\theta}\phi\right)  \mid_{A}=0$
and $\left|  \iota_{e_{1}\wedge...\wedge e_{r+1}}\left(  \phi|_{A}\right)
\right|  =1$ for any orthonormal tangent vector $e_{i}$'s on $A$.
\end{definition}

Remark: Recall that the volume of a calibrated submanifold is topological.\ In
this case, for the fundamental class $\left[  A\right]  \in H_{r+1}\left(
A,\mathbb{Z}\right)  $ and $\left[  Re\left(  e^{i\theta}\phi\right)  \right]
$ $\in H^{r+1}\left(  M,\mathbb{R}\right)  $,
\[
vol\left(  A\right)  =\int_{A}Re\left(  e^{i\theta}\phi\right)  =\left[
A\right]  \cdot\left[  Re\left(  e^{i\theta}\phi\right)  \right]
\]

Using the classification result of $\mathbb{C}$-VCPs in theorem
\ref{Classification of CVCP}, $\phi$ must be either a holomorphic volume form
in a Calabi-Yau manifold or a holomorphic symplectic form in a hyperk\"{a}hler
manifold. In the former case, instantons are called special Lagrangian
submanifolds (see e.g. \cite{Sch}). In the latter case,
\[
\operatorname{Re}\left(  e^{i\theta}\left(  \omega_{I}-i\omega_{K}\right)
\right)  =\cos\theta\omega_{I}+\sin\theta\omega_{K}=\omega_{\left(  \cos
\theta\right)  I+\left(  \sin\theta\right)  K}%
\]
so instantons are $J_{\theta}$-holomorphic curves where $J_{\theta}=\left(
\cos\theta\right)  I+\left(  \sin\theta\right)  K$ .

\bigskip

\subsection{Dirichlet and Neumann Branes}

Schoen's school studies (\cite{Butscher}, \cite{Qiu}, \cite{Sch}) free
boundary value problem for special Lagrangian submanifolds in Calabi-Yau
manifolds $M$. Suppose that $A$ is a special Lagrangian submanifold of zero
phase, i.e. calibrated by $\operatorname{Re}\Omega_{M}$, in $M$. If the
boundary of $A$ is non-empty and lies on a submanifold $C$ in $M$, then (i)
$C$ being a complex hypersurface in $M$ corresponds to Neumann boundary
condition on $A$ and (ii) $C$ being a special Lagrangian submanifold of phase
$\pi/2$ corresponds to Dirichlet boundary condition on $A$. Motivated from
these, we have the following definitions of branes in Hermitian complex manifolds.

\begin{definition}
On a Hermitian complex manifold $\left(  M,\omega\right)  $ of complex
dimension n with an r-fold $\mathbb{C}$-\textit{VCP} $\phi\in\Omega
^{r+1,0}\left(  M\right)  ,$

(i) a submanifold $C$ is called a $\mathbf{Neumann}$ $\mathbf{brane}$ (abbrev.
N-brane) if $\dim\left(  C\right)  $ $=$ $n+r-1$ and
\[
\phi|_{C}=0,
\]

(ii)\ an n-dimensional submanifold$\;C$\ is called a $\mathbf{Dirichlet}$
$\mathbf{brane}$ (abbrev. D-brane) with phase $\theta$\ $\in\mathbb{R}$ if
\[
\omega|_{C}=0\;,\;Re\left(  e^{i\theta}\phi\right)  \mid_{C}=0
\]
\end{definition}

Even though branes are defined in Hermitian complex manifolds, for simplicity
we will focus on branes in a closed K\"{a}hler manifold.

As in the real setting, N-branes are submanifolds in $M$ with the biggest
dimension on which $\phi$ vanishes. Furthermore N-branes are complex
submanifolds in $M$.

\begin{proposition}
Suppose that $M$ is an $n$-dimensional K\"{a}hler manifold with an $r$-fold
$\mathbb{C}$-VCP. Assume that $S$ is a submanifold in $M$ such that $\phi
|_{S}=0$, then $\dim S\leq n+r-1.$

When the equality sign holds, i.e. $S$ is a N-branes, then it is a complex
submanifold in $M$.
\end{proposition}

\begin{proof}
It suffices to check the linear case, i.e. $T_{x}S$, $x\in S$. Consider a
tangent vector $a$ $\in T_{x}S$ and $\tilde{a}=\left(  a-\sqrt{-1}Ja\right)
/\sqrt{2}$. Since$\;\phi\in\Omega^{r+1,0}$, $\iota_{\tilde{a}}\phi=\sqrt
{2}\iota_{a}\phi$.Therefore, for $b_{i}$'s in $T_{x}S.$
\[
\sqrt{2}\phi\left(  a,b_{1},...,b_{r}\right)  =\phi\left(  \tilde{a}%
,b_{1},...,b_{r}\right)  .
\]
So if $\phi\left(  a,b_{1},...,b_{r}\right)  =0$, then $\phi\left(
Ja,b_{1},...,b_{r}\right)  =0$, namely, $\phi|_{S}=0$ implies $\phi|_{S+JS}%
=0$. And more if $S$ has the maximal dimension, then $T_{x}S$ is complex
linear subspace. In that case, as in the real setting, we have $\dim
_{\mathbb{C}}S=\left(  n+r-1\right)  /2$.
\end{proof}

The above proposition can also be verified case by case using the
classification result of $\mathbb{C}$-VCP. The following classification of
instantons and branes for $\mathbb{C}$-\textit{VCP} is presented with the
classification of $\mathbb{C}$-\textit{VCP} on a closed K\"{a}hler manifold.

\bigskip

\textbf{Instantons and branes in CY manifold} \textbf{:} \textbf{(n-1)-fold
}$\mathbb{C}$\textbf{-}\textit{VCP}

A closed K\"{a}hler manifold of complex dimension n with $(n-1)$-fold
$\mathbb{C}$\textbf{-}\textit{VCP} $\phi$ is a $CY$ n-fold.

An instanton in the Calabi-Yau $n$-fold is a special Lagrangian submanifold
with phase $\theta\;$since it is calibrated by $\operatorname{Re}\left(
e^{i\theta}\phi\right)  $. As in the previous proposition, an $N$-brane in the
Calabi-Yau $n$-fold is a complex hypersurface, and a $D$-brane is a special
Lagrangian submanifold with phase $\theta-\pi/2,$ because it is calibrated by
$\operatorname{Re}\left(  e^{i\left(  \theta-\pi/2\right)  }\phi\right)  =$
$\operatorname{Im}$ $\left(  e^{i\theta}\phi\right)  $.

\bigskip

\textbf{Instantons and branes in hyperk\"{a}hler manifold :} \textbf{1-fold
}$C$\textbf{-}\textit{VCP}

A closed K\"{a}hler manifold of complex dimension $2n$ with a K\"{a}hler form
$\omega_{J}$ and a $1$-fold $C$\textbf{-}\textit{VCP} $\phi$ is a
hyperk\"{a}hler manifold. Denote $\phi=:\omega_{I}$\ $-\sqrt{-1}\omega_{K} $
and $J$ , $I$ and $K$ as the complex structure corresponding to K\"{a}hler
structures $\omega_{J}$, $\omega_{I}$ and $\omega_{K}$ , respectively. And
more, by putting $e^{i\theta}\phi$ in place of $\phi$, one can observe
$\operatorname{Re}e^{i\theta}\phi$ is another K\"{a}hler structure with a
complex structure $J_{\theta}=:$ $\cos\theta I+\sin\theta K$.

Now, an instanton in a hyperk\"{a}hler manifold is a $J_{\theta}$-holomorphic
curve since it is calibrated by $\operatorname{Re}e^{i\theta}\phi$, namely
preserved by $J_{\theta}.$ An $N$-brane in a hyperk\"{a}hler manifold is a
real $2n$-dimensional submanifold where $\phi$ vanishes, and as in the
previous proposition, it is equivalently a $J$-complex Lagrangian which is a
complex submanifold preserved by a complex structure $J$ with complex
dimension $n$. A $D$-brane is a real $2n$-dimensional submanifold where
$\omega$ and $\operatorname{Re}\left(  e^{i\theta}\phi\right)  $ vanish. One
can show that a $D$-brane is preserved by $J_{\theta+\pi/2}$ the almost
complex structure corresponding to $-\operatorname{Im}\left(  e^{i\theta}%
\phi\right)  $, i.e. $e^{i\theta}\phi=\omega_{J_{\theta}}$\ $-\sqrt{-1}%
\omega_{J_{\theta+\pi/2}}$. So a $D$-brane is equivalently, a $J_{\theta
+\pi/2}$-complex Lagrangian.

The above classification of instantons, $N$-branes and $D$-branes in manifolds
with $\mathbb{C}$-VCP is summarized in the table in page \pageref{Table C-VCP}.

\section{Symplectic Geometry on Knot Spaces}

Recall that a $1$-fold VCP on $M$ is a symplectic structure. In general an
$r$-fold VCP form on $M$ induces a symplectic structure on the space of
embedded submanifolds $\Sigma$ in $M$ of dimension $r-1$, which we simply call
a (multi-dimensional) knot space $\mathcal{K}_{\Sigma}M=Map\left(
\Sigma,M\right)  /Diff\left(  \Sigma\right)  $. For instance when $M$ is an
oriented three manifold, $\mathcal{K}_{\Sigma}M$ is the space of knots in $M$.
In this case Brylinski \cite{Bry} showed that $\mathcal{K}_{\Sigma}M$ has a
natural complex structure which makes it an infinite dimensional K\"{a}hler
manifold and used it to study the problem of geometric quantization. For
general $M$ with a VCP, we will relate the symplectic geometry of
$\mathcal{K}_{\Sigma}M$ to the geometry of branes and instantons in $M$.

When $M$ is a K\"{a}hler manifold with a $\mathbb{C}$-VCP, say a Calabi-Yau
manifold, one might try to complexify the above construction to define a
holomorphic symplectic structure (i.e. $1$-fold $\mathbb{C}$-VCP) on the
symplectic quotient $Map\left(  \Sigma,M\right)  //Diff\left(  \Sigma\right)
$. However there are various difficulties due to the fact that $Map\left(
\Sigma,M\right)  $ does not have a natural symplectic structure unless we fix
a background volume form on $\Sigma$. We will resolve this problem and define
a natural holomorphic symplectic structure on the isotropic knot space
$\mathcal{\hat{K}}_{\Sigma}M$. We do not know whether it is a hyperk\"{a}hler
structure or not because the symplectic form on $\mathcal{\hat{K}}_{\Sigma}M$
induced from the K\"{a}hler form on $M$ may not be closed.

We will also show that both complex hypersurfaces (i.e. N-branes) and special
Lagrangian submanifolds of phase $-\pi/2$ (i.e. D-branes) in $M$ correspond to
complex Lagrangian submanifolds in $\mathcal{\hat{K}}_{\Sigma}M$, but with
respect to different almost complex structures in the twistor $S^{2}$-family
for $\mathcal{\hat{K}}_{\Sigma}M$. For instance when $M$ is a three
dimensional Calabi-Yau manifold, $\mathcal{\hat{K}}_{\Sigma}M$ is roughly the
space of loops, or strings, in $M$ with an equivalence relation generated by
deformations along complex directions (see remark \ref{CY3 remark} for details).

\subsection{\label{Symplectic on Loop}Symplectic Structure on Knot Spaces}

Let $\left(  M,g\right)  $ be an $n$-dimensional Riemannian manifold $M$ with
a closed \textit{VCP} form $\phi$ of degree $r+1$. Suppose $\Sigma$ is any
$\left(  r-1\right)  $-dimensional oriented closed manifold $\Sigma$, we
consider the mapping space of embeddings from $\Sigma$ to $M$,
\[
Map\left(  \Sigma,M\right)  =\{f:\Sigma\rightarrow M\mid
f\;\text{is\ an\ embedding.}\}\text{.}%
\]
Let
\[
ev:\Sigma\times Map\left(  \Sigma,M\right)  \rightarrow M\;
\]
be the evaluation map $ev\left(  x,f\right)  =f(x)$ and $p_{1}$, $p_{2}$ be
the projection map from $\Sigma\times Map\left(  \Sigma,M\right)  $ to its
first and second factor respectively. We define a two form $\mathbb{\omega
}_{Map}$ on $Map\left(  \Sigma,M\right)  $ by taking the transgression of the
VCP form $\phi$,
\[
\omega_{Map}=\left(  p_{2}\right)  _{\ast}\left(  ev\right)  ^{\ast}\phi
=\int_{\Sigma}ev^{\ast}\phi\text{.}%
\]
Explicitly, suppose $u$ and $v$ are tangent vectors to $Map\left(
\Sigma,M\right)  $ at $f$, that is $u,v\in\Gamma\left(  \Sigma,f^{\ast}\left(
T_{M}\right)  \right)  $, we have
\[
\omega_{Map}\left(  u,v\right)  =\int_{\Sigma}\iota_{u\wedge v}\phi\text{.}%
\]

Since $\left(  \iota_{u\wedge v}\phi\right)  |_{\Sigma}$ can never be a top
degree form if $u$ is tangent to $\Sigma$, $\omega_{Map}$ degenerates along
tangent directions to the orbits of the natural action of $Diff\left(
\Sigma\right)  $ on $Map\left(  \Sigma,M\right)  $. Thus it descends to a two
form $\mathbf{\omega}^{\mathcal{K}}$ on the quotient space
\[
\mathcal{K}_{\Sigma}M=Map\left(  \Sigma,M\right)  /Diff\left(  \Sigma\right)
\text{,}%
\]
the space of submanifolds in $M$. For simplicity, we call it a
(multi-dimensional) \emph{knot space}. Note that tangent vectors to
$\mathcal{K}_{\Sigma}M$ are sections of the normal bundle of $\Sigma$ in $M$.

On $\mathcal{K}_{\Sigma}M$ there is a natural $L^{2}$-metric given as follow:
Suppose $\left[  f\right]  \in\mathcal{K}_{\Sigma}M$ and $u,v$ $\in
\Gamma(\Sigma,N_{\Sigma/M})$ are tangent vectors to $\mathcal{K}_{\Sigma}M$ at
$\left[  f\right]  $, then
\[
\mathbf{g}_{\left[  f\right]  }^{\mathcal{K}}\left(  u,v\right)  =%
{\displaystyle\int_{\Sigma}}
g\left(  u,v\right)  \nu_{\Sigma},
\]
where $\nu_{\Sigma}$ is the volume form of $\Sigma$ with respect to the
induced metric on $\Sigma$. We define an endomorphism $\mathbf{J}%
^{\mathcal{K}}$ on the tangent space of $\mathcal{K}_{\Sigma}M$ as follow:
Suppose $\left[  f\right]  \in\mathcal{K}_{\Sigma}M$ and $u,v$ $\in
\Gamma(\Sigma,N_{\Sigma/M})$ are tangent vectors to $\mathcal{K}_{\Sigma}M$ at
$\left[  f\right]  $, then
\[
\mathbf{\omega}_{\left[  f\right]  }^{\mathcal{K}}\left(  u,v\right)
=\mathbf{g}_{\left[  f\right]  }^{\mathcal{K}}\left(  \mathbf{J}_{\left[
f\right]  }^{\mathcal{K}}(u),v\right)  \text{.}%
\]

\begin{proposition}
Suppose $\left(  M,g\right)  $ is a Riemannian manifold with a VCP. Then
$\mathbf{J}^{\mathcal{K}}$ is a Hermitian almost complex structure on
$\mathcal{K}_{\Sigma}M$, i.e. a $1$-fold VCP.
\end{proposition}

\begin{proof}
Suppose $\left[  f\right]  \in\mathcal{K}_{\Sigma}M$ and $u,w$ $\in
\Gamma(\Sigma,N_{\Sigma/M})$ are tangent vectors to $\mathcal{K}_{\Sigma}M$ at
$\left[  f\right]  $. Let $e_{i}$'s be an oriented orthonormal base of
$\Sigma$ at a point $x$. Then
\begin{align*}
\imath_{u\wedge w}\left(  ev^{\ast}\phi\right)   &  =\phi\left(
e_{1},...e_{r-1},u,w\right)  v_{\Sigma}\\
&  =g\left(  \chi(e_{1},...e_{r-1},u),w\right)  v_{\Sigma}.
\end{align*}
From the relationship $\mathbf{\omega}_{\left[  f\right]  }^{\mathcal{K}%
}\left(  u,v\right)  =\mathbf{g}_{\left[  f\right]  }^{\mathcal{K}}\left(
\mathbf{J}_{\left[  f\right]  }^{\mathcal{K}}(u),w\right)  $ and the
definitions of $\mathbf{\omega}_{\left[  f\right]  }^{\mathcal{K}}$ and
$\mathbf{g}_{\left[  f\right]  }^{\mathcal{K}}$, we conclude that
$\mathbf{J}^{\mathcal{K}}$ on the tangent space of $\mathcal{K}_{\Sigma}M$ at
$\left[  f\right]  $ is given by
\[%
\begin{tabular}
[c]{llll}%
$\mathbf{J}^{\mathcal{K}}:$ & $\Gamma(\Sigma,N_{\Sigma/M})$ & $\rightarrow$ &
$\Gamma(\Sigma,N_{\Sigma/M})$\\
& $\mathbf{J}^{\mathcal{K}}\left(  u\right)  $ & $=$ & $\chi(e_{1}%
,...e_{r-1},u).$%
\end{tabular}
\]

On the other hand, from the remark \ref{Dimension of branes} on the dimension
of branes, we know that if $\chi$ is an $r$-fold VCP on $T_{x}M$, then
$\chi\left(  e_{1},...,e_{r-1,}\cdot\right)  $ defines a $1$-fold VCP on the
orthogonal complement to any oriented orthonormal vectors $e_{i}$'s at
$T_{x}M$. This implies that $\mathbf{J}^{\mathcal{K}}$ is a $1$-fold VCP on
$\mathcal{K}_{\Sigma}M$.
\end{proof}

\bigskip

Remark: The above proof actually shows that if $\mathbf{J}^{\mathcal{K}}$ on
$\mathcal{K}_{\Sigma}M$ is induced by an arbitrary differential form $\phi$ on
$M$, then $\mathbf{J}^{\mathcal{K}}$ is an Hermitian almost complex structure
on $\mathcal{K}_{\Sigma}M$ if and only if $\phi$ is a $r$-fold VCP form on $M$.

\bigskip

When the VCP form $\phi$ on $M$ is closed and $\Sigma$ is a closed manifold,
then $\mathbb{\omega}_{Map}$ on $Map(\Sigma,M)$ is also closed because
\[
d\mathbb{\omega}_{Map}\mathbf{=}d%
{\displaystyle\int_{\Sigma}}
ev^{\ast}\phi=%
{\displaystyle\int_{\Sigma}}
ev^{\ast}d\phi\;=0.
\]
Therefore $\omega_{Map}$ descends to a closed $1$-fold VCP form on
$\mathcal{K}_{\Sigma}M$.

Conversely the closedness of $\omega_{Map}$ on $Map\left(  \Sigma,M\right)  $,
or on $\mathcal{K}_{\Sigma}M$, implies the closedness of $\phi$ on $M.$ Since
this type of localization arguments will be used several times in this paper,
we include the proof of the following standard lemma.

\begin{lemma}
\label{Localize Lemma}(\textbf{Localization}) Let $\Sigma$ be an
$s$-dimensional manifold without boundary. A$\;$form $\eta$ of degree $k>s$ on
a manifold $M$ vanishes if the corresponding $\left(  k-s\right)  $-form on
$Map\left(  \Sigma,M\right)  \;$obtained by transgression vanishes.
\end{lemma}

\begin{proof}
We need to show that
\[
\eta(p)(v_{1},v_{2},..,v_{k})=0
\]
for any fixed $p$ $\in$ $M$ and any fixed $v_{i}\in T_{p}M$. For simplicity,
we may choose $v_{i}$'s to be orthonormal vectors. We can find $f\in$
$Map(\Sigma,M)$ with $p\in f\left(  \Sigma\right)  $ such that $v_{1}%
,...,v_{s}$ are along the tangential directions, and $v_{s+1},...,v_{k}$ are
along the normal directions at $p$ in $f(\Sigma)$. Moreover, we can choose
sections $\check{v}_{s+1},\tilde{v}_{s+2},...,\tilde{v}_{k}$ $\in$
$\Gamma(\Sigma,f^{\ast}(T_{M}))$ which equal $v_{s+1},...,v_{k}$ at $p$
respectively. By multipling $\check{v}_{s+1}$ with a sequence of functions on
$\Sigma$ approaching the delta function at $p$, we obtain sections $\left(
\check{v}_{s+1}\right)  _{\varepsilon}$ which approach $\delta\left(
p\right)  v_{s+1}$ as $\varepsilon\rightarrow0$ where $\delta\left(  p\right)
$\ is Dirac delta function. Therefore,
\[
\eta(p)(v_{1},v_{2},..,v_{k})=\underset{\varepsilon\rightarrow0}{lim}(%
{\displaystyle\int_{\Sigma}}
ev^{\ast}\eta)(\left(  \check{v}_{s+1}\right)  _{\varepsilon},\tilde{v}%
_{s+2},..,\tilde{v}_{k}).
\]
From the given condition $%
{\displaystyle\int_{\Sigma}}
ev^{\ast}\eta=0$, we conclude that
\[
\eta(p)(v_{1},v_{2},..,v_{k})=0.
\]
Hence the result.
\end{proof}

As a corollary of the above lemma and discussions, we have the following result.

\begin{proposition}
Suppose $\left(  M,g\right)  $ is a Riemannian manifold with a differential
form $\phi$ of degree $r+1$. Then $\phi$ is a closed $r$-fold VCP form on $M$
if and only if $\mathbf{\omega}^{\mathcal{K}}$ is an almost K\"{a}hler
structure on $\mathcal{K}_{\Sigma}M$, i.e. a closed $1$-fold VCP form.
\end{proposition}

Remark: In general the transgression of any closed form $\phi$ always gives a
closed form on the knot space $\mathcal{K}_{\Sigma}M$ of degree $r+1-s$ where
$\dim\Sigma=s$. However this can never be a VCP form unless $s=r-1$. To see
this, we can use the above localization method and the fact that there is no
VCP form of degree bigger than two on any vector space with sufficiently large dimension.

When $r=1$, that is $M$ is a symplectic manifold, $Map\left(  \Sigma,M\right)
$ is the same as $M$ and therefore it is symplectic by trivial reasons.

\bigskip

Remark: Given any $\phi$-preserving vector field $v\in Vect\left(
M,\phi\right)  $, it induces a vector field $V$ of $\mathcal{K}_{\Sigma}M$
preserving $\mathbf{\omega}^{\mathcal{K}}$. Furthermore, if $v$ is $\phi
$-Hamiltonian, that is
\[
\iota_{v}\phi=d\eta\text{,}%
\]
for some $\eta\in\Omega^{r}\left(  M\right)  $, then we can define a function
on the knot space,
\begin{align*}
F_{\eta}  &  :\mathcal{K}_{\Sigma}M\rightarrow\mathbb{R}\\
F_{\eta}\left(  f\right)   &  =\int_{\Sigma}f^{\ast}\eta\text{.}%
\end{align*}
It is easy to see that $F_{\eta}$ is a Hamiltonian function on the symplectic
manifold $\mathcal{K}_{\Sigma}M$ whose Hamiltonian vector field equals $V$.

\subsection{Holomorphic Curves and Lagrangians in Knot Spaces}

In this subsection, we are going to show that holomorphic disks (resp.
Lagrangian submanifolds) in the knot space $\mathcal{K}_{\Sigma}M$ of $M$
correspond to instantons (resp. branes) in $M$. More generally the geometry of
vector cross products on $M$ should be closely related to the symplectic
geometry of its knot space $\mathcal{K}_{\Sigma}M$. A natural problem is to
understand the analog of the Floer's Lagrangian intersection theory for
manifolds with vector cross products. Note that if two branes $C_{1}$ and
$C_{2}$ intersect transversely along a submanifold $\Sigma$, then the
dimension of $\Sigma$ equals $r-1$ and $\left[  \Sigma\right]  $ represents a
transverse intersection point of (Lagrangians) $\mathcal{K}_{\Sigma}C_{1}$ and
$\mathcal{K}_{\Sigma}C_{2}$ in $\mathcal{K}_{\Sigma}M$. The converse is also true.

Because $\mathcal{K}_{\Sigma}M$ is of infinite dimensional, a Lagrangian
submanifold is defined as a subspace in $\mathcal{K}_{\Sigma}M$ where the
restriction of $\mathbf{\omega}^{\mathcal{K}}$ vanishes and with the property
that any vector field $\mathbf{\omega}^{\mathcal{K}}$-orthogonal to $L$ is a
tangent vector field along $L$, see \cite{M} and \cite{Bry}. In fact, it is
easy to see that the latter condition for a submanifold $L$ to be Lagrangian
is equivalent to the statement that $\mathbf{\omega}^{\mathcal{K}}$ will not
vanish on any bigger space containing $L$. So, we refer the this condition as
the \textbf{maximally} \textbf{self }$\mathbf{\omega}^{\mathcal{K}}%
$\textbf{-perpendicular} condition.

\begin{proposition}
\label{ProD}Suppose that $M$ is an $n$-dimensional manifold $M$ with a closed
\textit{VCP} form $\phi$ of degree $r+1$ and $C$ is a submanifold in $M$. Then
$\mathcal{K}_{\Sigma}C$ is a Lagrangian submanifold in $\mathcal{K}_{\Sigma
}M\;$if and only if $C$ is a brane in $M$.
\end{proposition}

\begin{proof}
For the \textit{only if} part, we suppose that $\mathcal{K}_{\Sigma}C\;$is a
Lagrangian in $\mathcal{K}_{\Sigma}M$ and we want to show that $\phi$ vanishes
along $C$ and dim $C=\left(  n+r-1\right)  /2$. For any fixed \ $\left[
f\right]  $ $\in$ $\mathcal{K}_{\Sigma}C$, we have
\[
0=\mathbf{\omega}_{\left[  f\right]  }^{\mathcal{K}}(u,v)=%
{\displaystyle\int_{\Sigma}}
\iota_{u\wedge v}\left(  ev^{\ast}\phi\right)  \;
\]
where $u,v$ $\in\Gamma(\Sigma,N_{\Sigma/C})$. By applying the localization
method as in lemma \ref{Localize Lemma} along $C$, we obtain the following:
For any $x$ in $\Sigma$,
\[
\phi(u\left(  x\right)  ,v\left(  x\right)  ,e_{1},...,e_{r-1})=0,
\]
where $e_{1,}...,e_{r-1\;}$ are any oriented orthonormal vectors of
$T_{x}\Sigma$. By varying $u,v,$ $f$ to cover $T_{C}.$ One can show $\phi
\mid_{C}=0$.

Moreover, $\mathcal{K}_{\Sigma}C$ being $\mathbf{\omega}^{\mathcal{K}}%
$-perpendicular in $\mathcal{K}_{\Sigma}M$ implies that $C$ has the biggest
possible dimension with $\phi|_{C}=0$. As explained in the remark following
the definition of branes that this gives $\dim$ $C=\left(  n+r-1\right)  /2$.

For the\textit{ if} part, we assume that $C$ is any brane in $M$. The
condition $\phi|_{C}=0$ implies that $\mathbf{\omega}^{\mathcal{K}%
}|_{\mathcal{K}_{\Sigma}C}=0$. To show that $\mathcal{K}_{\Sigma}C\;$is a
Lagrangian in $\mathcal{K}_{\Sigma}M$, we need to verify the maximally self
$\mathbf{\omega}^{\mathcal{K}}$-perpendicular condition in $\mathcal{K}%
_{\Sigma}M$. Recall that the tangent space at any point $\left[  f\right]
\in\mathcal{K}_{\Sigma}C$ is $\Gamma(\Sigma,N_{\Sigma/C})$. Suppose there is a
section $v$ in $\Gamma(\Sigma,N_{\Sigma/M})$ but not in $\Gamma(\Sigma
,N_{\Sigma/C})$ such that it is $\mathbf{\omega}^{\mathcal{K}}$-perpendicular
to $\Gamma(\Sigma,N_{\Sigma/C})$. By the localization arguments as in lemma
\ref{Localize Lemma}, given any point $x\in f\left(  \Sigma\right)  \subset
C$, $\phi$ vanishes on the linear space spanned by $v\left(  x\right)  $ and
$T_{x}C$. This contradicts to the fact that $C$ is a brane. Hence the result.
\end{proof}

Remark: $\phi$-Hamiltonian deformation of a brane $C$ in $M$ corresponds to
Hamiltonian deformation of the corresponding Lagrangian submanifold
$\mathcal{K}_{\Sigma}C$ in the symplectic manifold $\mathcal{K}_{\Sigma}M$.
More precisely, if $v$ is a normal vector field to $C$ satisfying
\[
\iota_{v}\phi=d\eta,
\]
for some $\eta\in\Omega^{r-1}\left(  C\right)  $, then the transgression of
$\eta$ defines a function on $\mathcal{K}_{\Sigma}C$ which generates a
Hamiltonian deformation of $\mathcal{K}_{\Sigma}C$.

Next we discuss holomorphic disks, i.e. instantons, in $\mathcal{K}_{\Sigma
}M.$ We consider a two dimensional disk $D$ in $Map\left(  \Sigma,M\right)  $
such that for each tangent vector $v\in T_{\left[  f\right]  }D$, the
corresponding vector field in $\Gamma\left(  \Sigma,f^{\ast}T_{M}\right)  $ is
normal to $\Sigma$. We call such a disk $D$ as a \textbf{normal disk. }For
simplicity we assume that the $r+1$ dimensional submanifold
\[
A=\underset{f\in D}{\;{\large \cup}}f\left(  \Sigma\right)  \subset M\text{,}%
\]
is an embedding. This is always the case if $D$ is small enough. Notice that
$A$ is diffeomorphic to $D\times\Sigma$. We will denote the corresponding disk
in $\mathcal{K}_{\Sigma}M$ as $\hat{D}=:\pi\left(  D\right)  $. We remark that
the principal fibration
\[
Diff\left(  \Sigma\right)  \rightarrow Map\left(  \Sigma,M\right)
\overset{\pi}{\rightarrow}\mathcal{K}_{\Sigma}M
\]
has a canonical connection (see \cite{Bry}) and $D$ being a normal disk is
equivalent to it being an integral submanifold for the horizontal distribution
of this connection.

In the following proposition, we describe the relation between a disk $\hat
{D}$ in $\mathcal{K}_{\Sigma}M$ given above and the corresponding $r+1$
dimensional subspace $A$ in $M$.

\begin{proposition}
\label{ProDD}Suppose that $M$ is a manifold with a closed $r$-fold VCP form
$\phi$ and $\mathcal{K}_{\Sigma}M$ is its knot space as before. For a normal
disk $D$ in $Map\left(  \Sigma,M\right)  $, $\hat{D}=:\pi\left(  D\right)  $
is a $J^{\mathcal{K}}$-holomorphic disk in $\mathcal{K}_{\Sigma}M$, i.e.
calibrated by $\omega^{\mathcal{K}}$, if and only if $A$ is an instanton in
$M$ and $A\rightarrow D$ is a Riemannian submersion.
\end{proposition}

We call such an $A$ a \emph{tight instanton.}

\begin{proof}
For a fixed $\left[  f\right]  \in\hat{D}$, we consider $\nu,\mu\in T_{\left[
f\right]  }(\hat{D})\subset\Gamma(\Sigma,N_{\Sigma/A})$. Since $\phi$ is a
calibrating form, we have,
\[
\phi(\nu,\mu,e_{1},...,e_{r-1})\leq Vol_{A}(\nu,\mu,e_{1},...,e_{r-1})=\left|
\nu\wedge\mu\right|
\]
where $e_{1},...,e_{r-1}$ is any orthonormal frame on $f\left(  \Sigma\right)
$. In particular we have%
\[%
{\displaystyle\int_{f\left(  \Sigma\right)  }}
\iota_{\nu\wedge\mu}\left(  ev^{\ast}\phi\right)  \leq%
{\displaystyle\int_{f\left(  \Sigma\right)  }}
\left|  \nu\wedge\mu\right|  vol_{\Sigma},
\]
and the equality sign holds for every $\left[  f\right]  \in\hat{D}$ if and
only if $A$ is an instanton in $M$. We will simply denote $%
{\displaystyle\int_{f\left(  \Sigma\right)  }}
$ by $%
{\displaystyle\int_{\Sigma}}
$. Notice that the symplectic form on $\mathcal{K}_{\Sigma}M$ is given by,
\[
\mathbf{\omega}_{\left[  f\right]  }^{\mathcal{K}}(\nu,\mu)=%
{\displaystyle\int_{\Sigma}}
\iota_{\nu\wedge\mu}\left(  ev^{\ast}\phi\right)  \text{.}%
\]

Since $\mathbf{\omega}^{\mathcal{K}}$ is a 1-fold VCP form on $\mathcal{K}%
_{\Sigma}M$, we have
\[
\mathbf{\omega}_{\left[  f\right]  }^{\mathcal{K}}(\nu,\mu)\leq\left(  \left|
\nu\right|  _{{\LARGE \mathcal{K}}}^{2}\left|  \mu\right|
_{{\LARGE \mathcal{K}}}^{2}-\left\langle \nu,\mu\right\rangle
_{{\LARGE \mathcal{K}}}^{2}\right)  ^{1/2}%
\]
where $\left\langle \nu,\mu\right\rangle _{\mathcal{K}}=:g_{\mathcal{K}%
}\left(  \nu,\mu\right)  $ and $\left|  a\right|  _{\mathcal{K}}%
^{2}=:g_{\mathcal{K}}\left(  a,a\right)  $. Furthermore the equality sign
holds when $\hat{D}\;$is a$\;J^{\mathcal{K}}$- holomorphic disk
in$\;\mathcal{K}_{\Sigma}M$.

To prove the \textit{only if }part, we suppose that $\hat{D}\;$is
a$\;J^{\mathcal{K}}$-holomorphic disk in$\;\mathcal{K}_{\Sigma}M$. From above
discussions, we have
\begin{align*}%
{\displaystyle\int_{\Sigma}}
\left|  \nu\wedge\mu\right|   &  \geq%
{\displaystyle\int_{\Sigma}}
\imath_{\nu\wedge\mu}\left(  ev^{\ast}\phi\right)  =\left(  \left|
\nu\right|  _{{\LARGE \mathcal{K}}}^{2}\left|  \mu\right|
_{{\LARGE \mathcal{K}}}^{2}-\left\langle \nu,\mu\right\rangle
_{{\LARGE \mathcal{K}}}^{2}\right)  ^{1/2}\\
&  =\left(
{\displaystyle\int_{\Sigma}}
\left|  \nu\right|  ^{2}%
{\displaystyle\int_{\Sigma}}
\left|  \mu\right|  ^{2}-\left(
{\displaystyle\int_{\Sigma}}
\left\langle \nu,\mu\right\rangle \right)  ^{2}\right)  ^{1/2}.
\end{align*}
Combining with the H\"{o}lder inequality,
\[
\left(
{\displaystyle\int_{\Sigma}}
\left|  \nu\wedge\mu\right|  \right)  ^{2}+\left(
{\displaystyle\int_{\Sigma}}
\left\langle \nu,\mu\right\rangle \right)  ^{2}\leq%
{\displaystyle\int_{\Sigma}}
\left|  \nu\right|  ^{2}%
{\displaystyle\int_{\Sigma}}
\left|  \mu\right|  ^{2},
\]
we obtain%
\[%
\begin{array}
[c]{ccc}%
\text{(i)} &
{\displaystyle\int_{\Sigma}}
\iota_{\nu\wedge\mu}\left(  ev^{\ast}\phi\right)  =%
{\displaystyle\int_{\Sigma}}
\left|  \nu\wedge\mu\right|  vol_{\Sigma} & \text{and}\\
\text{(ii)} & \left(
{\displaystyle\int_{\Sigma}}
\left|  \nu\wedge\mu\right|  \right)  ^{2}+\left(
{\displaystyle\int_{\Sigma}}
\left\langle \nu,\mu\right\rangle \right)  ^{2}=%
{\displaystyle\int_{\Sigma}}
\left|  \nu\right|  ^{2}%
{\displaystyle\int_{\Sigma}}
\left|  \mu\right|  ^{2}. &
\end{array}
\]

Condition (i) says that $A$ is an instanton in $M$. Condition (ii) implies
that given any $\left[  f\right]  $, there exists constant $C_{1}$ and $C_{2}$
such that for any $x\in\Sigma$,
\[
\left|  \nu\left(  x\right)  \right|  =C_{1}\left|  \mu\left(  x\right)
\right|  \;,\;\angle\left(  \nu\left(  x\right)  ,\mu\left(  x\right)
\right)  =C_{2}.
\]
This implies that $A\rightarrow D\;$is a$\;$\emph{Riemannian\ submersion.}

For the \textit{if }part, we notice that $A$ being an instanton in $M$ implies
that
\begin{align*}%
{\displaystyle\int_{\Sigma}}
\left|  \nu\wedge\mu\right|   &  =%
{\displaystyle\int_{\Sigma}}
\iota_{\nu\wedge\mu}\left(  ev^{\ast}\phi\right)  =\mathbf{\omega}_{\left[
f\right]  }^{\mathcal{K}}(\nu,\mu)\leq\left(  \left|  \nu\right|
_{{\LARGE \mathcal{K}}}^{2}\left|  \mu\right|  _{{\LARGE \mathcal{K}}}%
^{2}-\left\langle \nu,\mu\right\rangle _{{\LARGE \mathcal{K}}}^{2}\right)
^{1/2}\\
&  =\left(
{\displaystyle\int_{\Sigma}}
\left|  \nu\right|  ^{2}%
{\displaystyle\int_{\Sigma}}
\left|  \mu\right|  ^{2}-\left(
{\displaystyle\int_{\Sigma}}
\left\langle \nu,\mu\right\rangle \right)  ^{2}\right)  ^{1/2}.
\end{align*}
Recall that the Riemannian submersion condition implies an equality ,
\[
\left(
{\displaystyle\int_{\Sigma}}
\left|  \nu\wedge\mu\right|  \right)  ^{2}+\left(
{\displaystyle\int_{\Sigma}}
\left\langle \nu,\mu\right\rangle \right)  ^{2}=%
{\displaystyle\int_{\Sigma}}
\left|  \nu\right|  ^{2}%
{\displaystyle\int_{\Sigma}}
\left|  \mu\right|  ^{2}%
\]
so the above inequality turns into an equality so that it gives
\[
\mathbf{\omega}_{\left[  f\right]  }^{\mathcal{K}}(\nu,\mu)=%
{\displaystyle\int_{\Sigma}}
\iota_{\nu\wedge\mu}\left(  ev^{\ast}\phi\right)  =\left(  \left|  \nu\right|
_{{\LARGE \mathcal{K}}}^{2}\left|  \mu\right|  _{{\LARGE \mathcal{K}}}%
^{2}-\left\langle \nu,\mu\right\rangle _{{\LARGE \mathcal{K}}}^{2}\right)
^{1/2}%
\]
i.e. $\hat{D}\;$is $J^{\mathcal{K}}\;$holomorphic in$\;\mathcal{K}_{\Sigma}M$.
\end{proof}

\section{Isotropic Knot Spaces of CY manifolds}

Recall that any volume form on a manifold $M$ of dimension $n$ determines a
closed $1$-fold VCP on the knot space $\mathcal{K}_{\Sigma}M=Map\left(
\Sigma,M\right)  /Diff\left(  \Sigma\right)  $ where $\Sigma$ is any closed
manifold of dimension $n-2$. It is natural to guess that the holomorphic
volume form on any Calabi-Yau $n$-fold $M$ would determine a closed $1$-fold
$\mathbb{C}$-VCP on the symplectic quotient $Map\left(  \Sigma,M\right)
//Diff\left(  \Sigma\right)  $. However we need to choose a background volume
form $\sigma$ on $\Sigma$ to construct a symplectic structure in $Map\left(
\Sigma,M\right)  $ which is only invariant under $Diff\left(  \Sigma
,\sigma\right)  $, the group of volume preserving diffeomorphisms of $\Sigma$
(\cite{D}\cite{H1}). Therefore we can only construct the symplectic quotient
$Map\left(  \Sigma,M\right)  //Diff\left(  \Sigma,\sigma\right)  =\mu
^{-1}\left(  0\right)  /Diff\left(  \Sigma,\sigma\right)  $. This larger space
is not hyperk\"{a}hler because the holomorphic two form $\int_{\Sigma}%
\Omega_{M}$ degenerates. We will show that a further symplectic reduction will
produce a Hermitian integrable complex manifold $\mathcal{\hat{K}}_{\Sigma}M$
with a 1-fold $\mathbb{C}$-VCP, in particular a holomorphic symplectic
structure. This may not be hyperk\"{a}hler because even the Hermitian complex
structure is integrable on $\mathcal{\hat{K}}_{\Sigma}M$, its K\"{a}hler form
may not closed. $\mathcal{\hat{K}}_{\Sigma}M$ should be regarded as a modified
construction for the non-existing space $Map\left(  \Sigma,M\right)
//Diff\left(  \Sigma\right)  $. We call this an \textit{isotropic knot space
}of $M$ for reasons which will become clear later.

Furthermore we will relate instantons, $N$-branes and $D$-branes in a
Calabi-Yau manifold $M$ to holomorphic curves and complex Lagrangian
submanifolds in the holomorphic symplectic manifold $\mathcal{\hat{K}}%
_{\Sigma}M$. These constructions is of particular interest when $M$ is a
Calabi-Yau threefold (see remark \ref{CY3 remark} for details).

\subsection{\label{SecConstrModifLoop}Holomorphic Symplectic Structures on
Isotropic Knot Spaces}

Let $M$ be a Calabi-Yau $n$-fold with a holomorphic volume form $\Omega_{M}$,
i.e. a closed $(n-1)$-fold $\mathbb{C}$-\textit{VCP}. To obtain a $1$-fold
$\mathbb{C}$-\textit{VCP} on a certain knot space by transgression, we first
fix an $n-2$ dimensional manifold $\Sigma\;$without boundary and we let
$Map\left(  \Sigma,M\right)  $ be the space of embeddings from $\Sigma$ to $M$
as before. For simplicity we assume that the first Betti number of $\Sigma$ is
zero, $b_{1}\left(  \Sigma\right)  =0$.

If we fix a background volume form $\sigma$ on $\Sigma$, then the K\"{a}hler
form $\omega$ on $M$ induces a natural symplectic form on $Map\left(
\Sigma,M\right)  $ as follow: for any tangent vectors $X$ and $Y$ at a point
$f$ in $Map\left(  \Sigma,M\right)  $, we define
\[
\mathbb{\omega}_{Map}\left(  X,Y\right)  =%
{\displaystyle\int_{\Sigma}}
\omega{\large (}X,Y{\large )}\sigma
\]
where$\,\,X,Y\in\Gamma(\Sigma,f^{\ast}(TM))$. Note that this symplectic
structure on $Map\left(  \Sigma,M\right)  $ is not invariant under general
diffeomorphisms of $\Sigma$. Instead it is preserved by the natural action of
$Diff\left(  \Sigma,\sigma\right)  $, the group of volume preserving
diffeomorphisms on $\left(  \Sigma,\sigma\right)  $.

As studied by Donaldson in \cite{D} and Hitchin in \cite{H1}, this action is
Hamiltonian on the components of $Map\left(  \Sigma,M\right)  $ consisting of
those $f$'s satisfying
\[
f^{\ast}\left(  \left[  \omega\right]  \right)  =0\in H^{2}\left(
\Sigma,\mathbb{R}\right)  ,
\]
and the moment map is given by
\[%
\begin{tabular}
[c]{llll}%
$\mu:$ & $Map_{0}\left(  \Sigma,M\right)  $ & $\rightarrow$ & $\Omega
^{1}\left(  \Sigma\right)  \diagup d\Omega^{0}\left(  \Sigma\right)  $\\
& $\;\;\;\;\;\;f$ & $\mapsto$ & $\;\;\;\;\;\mu\left(  f\right)  =\alpha\;$%
\end{tabular}
\]
for any one form $\alpha\in\Omega^{1}\left(  \Sigma\right)  $ satisfying
$d\alpha=f^{\ast}\omega$. Note that the dual of the Lie algebra of
$Diff\left(  \Sigma,\sigma\right)  $ can be naturally identified with
$\Omega^{1}\left(  \Sigma\right)  \diagup d\Omega^{0}\left(  \Sigma\right)  $.
In particular $\mu^{-1}\left(  0\right)  $ consists of isotropic embeddings of
$\Sigma$ in $M$. Therefore the symplectic quotient
\[
Map\left(  \Sigma,M\right)  //Diff\left(  \Sigma,\sigma\right)  =\mu
^{-1}\left(  0\right)  /Diff\left(  \Sigma,\sigma\right)
\]
is almost the same as the moduli space of isotropic submanifolds in $M$, which
is $\mu^{-1}\left(  0\right)  /Diff\left(  \Sigma\right)  $. Observe that the
definition of the map $\mu$ is independent of the choice of $\sigma$ and
therefore $\mu^{-1}\left(  0\right)  $ is preserved by the action of
$Diff\left(  \Sigma\right)  $.

Remark: For the moment map $\mu$ to be well-defined, the condition
$b_{1}\left(  \Sigma\right)  =0$ is necessary. However even in the case of
$b_{1}\left(  \Sigma\right)  \neq0$, which always happens when $M$ is a
Calabi-Yau threefold, there is modification of the symplectic quotient
construction and we can obtain a symplectic manifold which is a torus bundle
over the moduli space of isotropic submanifolds in $M$ with fiber dimension
$b_{1}\left(  \Sigma\right)  $ (\cite{D}\cite{H1}).

In order to construct a holomorphic symplectic manifold, we first consider a
complex closed $2$-form $\Omega_{Map}$ on $Map\left(  \Sigma,M\right)  $
induced from the holomorphic volume form $\Omega_{M}$ on $M$\textit{\ }by
transgression.
\begin{align*}
\Omega_{Map}  &  =\int_{\Sigma}ev^{\ast}\Omega_{M}\\
&  =\int_{\Sigma}ev^{\ast}\operatorname{Re}\Omega_{M}+\sqrt{-1}\int_{\Sigma
}ev^{\ast}\operatorname{Im}\Omega_{M}\\
&  =\omega_{I}-\sqrt{-1}\omega_{K}\text{,}%
\end{align*}
where $\omega_{I}=\int_{\Sigma}ev^{\ast}\,\left(  \operatorname{Re}\Omega
_{M}\right)  $ and $\omega_{K}=-\int_{\Sigma}ev^{\ast}\,\left(
\operatorname{Im}\Omega_{M}\right)  $. We also define endomorphisms $I\;$and
$K$ on $T_{f}\left(  \mu^{-1}\left(  0\right)  \right)  $ as follows:
\[
\omega_{I}\left(  A,B\right)  =g(IA,B),\;\omega_{K}\left(  A,B\right)
=g(KA,B)\;
\]
where $A,B\in T_{f}\left(  \mu^{-1}\left(  0\right)  \right)  $ and $g$ is the
natural $L^{2}$-metric on $\mu^{-1}\left(  0\right)  $. Both $\omega_{I}$ and
$\omega_{K}$ are degenerated two forms. We will define the isotropic knot
space as the symplectic reduction of $\mu^{-1}\left(  0\right)  $ with respect
to either $\omega_{I}$ or $\omega_{K}$ and show that it has a natural
holomorphic symplectic structure.

To do that, we define a distribution $\mathcal{D}$ on $\mu^{-1}\left(
0\right)  $ by
\[
\mathcal{D}_{f}=\left\{  X\in T_{f}\left(  \mu^{-1}\left(  0\right)  \right)
\subset\Gamma\left(  \Sigma,f^{\ast}\left(  T_{M}\right)  \right)  \mid
\iota_{X}\omega_{I,f}=0\right\}  \subset T_{f}\left(  \mu^{-1}\left(
0\right)  \right)  ,
\]
for any $f\in\mu^{-1}\left(  0\right)  $.

\begin{lemma}
For any $f\in\mu^{-1}\left(  0\right)  \subset Map\left(  \Sigma,M\right)  ,$
we have,
\[
\mathcal{D}_{f}=\Gamma\left(  \Sigma,T_{\Sigma}+JT_{\Sigma}\right)  .
\]
\end{lemma}

\begin{proof}
First, it is easy to see that $\mathcal{D}_{f}\supset\Gamma\left(
\Sigma,T_{\Sigma}\right)  $ since for any $X$ in $\Gamma\left(  \Sigma
,T_{\Sigma}\right)  $ and any $Y\;$in$\;\Gamma\left(  \Sigma,T_{M}\right)  $,
$\iota_{X\wedge Y}\operatorname{Re}\,\Omega_{M}$ can not be a top degree form
on $\Sigma$. By similar reasons, $\mathcal{D}_{f}$ is preserved by the
hermitian complex structure $J\;$on $M.$ Because $f$ is isotropic, we have
$\mathcal{D}_{f}\supset\Gamma\left(  \Sigma,T_{\Sigma}+JT_{\Sigma}\right)  $.
Now, as in lemma \ref{Localize Lemma}, we consider localization of
\[
0=\iota_{X}\omega_{I,f}=\int_{\Sigma}\iota_{X}ev^{\ast}\operatorname{Re}%
\,\Omega_{M}\;,
\]
at $x$ in $\Sigma$ by varying $\Sigma$ and we obtain%
\[
0=\iota_{X\left(  x\right)  \wedge E_{1}\wedge...\wedge E_{n-2}}%
\operatorname{Re}\,\Omega_{M}%
\]
where $E_{1,}...,E_{n-2}$ is an orthonormal oriented basis of $\left(
T_{\Sigma}\right)  _{x}$. This implies that
\[
\mathcal{D}_{f}=\Gamma\left(  \Sigma,T_{\Sigma}+JT_{\Sigma}\right)  ,
\]
because for any $X$ in $T_{M}\backslash\left(  T_{\Sigma}+JT_{\Sigma}\right)
$, there is $x$ in $\Sigma$ such that
\[
\iota_{X\left(  x\right)  \wedge E_{1}\wedge...\wedge E_{n-2}}%
\operatorname{Re}\,\Omega_{M}\neq0.
\]
Note that the same construction applied to $\omega_{K}$, instead of
$\omega_{I}$, will give another distribution which is identical to
$\mathcal{D}_{f}$ because $f$ is isotropic.
\end{proof}

Observe that, for each $f\in$ $\mu^{-1}\left(  0\right)  $,$\;$the rank of the
subbundle $T_{\Sigma}+JT_{\Sigma}$ in $f^{\ast}T_{M}$ is $2\left(  n-2\right)
$, i.e. constant rank on $\mu^{-1}\left(  0\right)  $. That is $\omega_{I}$ is
a closed two form of constant rank on $\mu^{-1}\left(  0\right)  $. From the
standard theory of symplectic reduction, the distribution $\mathcal{D}$ is
integrable and the space of leaves has a natural symplectic form descended
from $\omega_{I}$. We call this space as the \textbf{isotropic knot space} of
$M$ and we denote it as
\[
\mathcal{\hat{K}}_{\Sigma}M=\mu^{-1}\left(  0\right)  /\left\langle
\mathcal{D}\right\rangle
\]
where $\left\langle \mathcal{D}\right\rangle $ are equivalence relations
generated by the distribution $\mathcal{D}$.

\bigskip

Remark: The isotropic knot space $\mathcal{\hat{K}}_{\Sigma}M$ is a quotient
space of $\mu^{-1}\left(  0\right)  /Diff\left(  \Sigma\right)  $, the space
of isotropic submanifolds in $M$. If $f:\Sigma\rightarrow M$ parametrizes an
isotropic submanifold in $M$, then deforming $\Sigma$ along $T_{\Sigma}$
directions simply changes the parametrization of the submanifold $f\left(
\Sigma\right)  $, namely the equivalence class in $\mu^{-1}\left(  0\right)
/Diff\left(  \Sigma\right)  $ remains unchanged. However if we deform $\Sigma$
along $T_{\Sigma}\otimes\mathbb{C}$ directions, then their equivalence classes
in $\mathcal{\hat{K}}_{\Sigma}M$ remains constant. Notice that the isotropic
condition on $\Sigma$ implies that $T_{\Sigma}\otimes\mathbb{C}\cong
T_{\Sigma}\oplus J\left(  T_{\Sigma}\right)  \subset f^{\ast}T_{M}$.
Therefore, roughly speaking, $\mathcal{\hat{K}}_{\Sigma}M$ is the space of
isotropic submanifolds in $M$ divided by $Diff\left(  \Sigma\right)
\otimes\mathbb{C}$. In particular, the tangent space of $\mathcal{\hat{K}%
}_{\Sigma}M$ is given by
\[
T_{\left[  f\right]  }\left(  \mathcal{\hat{K}}_{\Sigma}M\right)
=\Gamma\left(  \Sigma,f^{\ast}\left(  TM\right)  /\left(  T_{\Sigma
}+JT_{\Sigma}\right)  \right)
\]
for any $\left[  f\right]  \;$in $\mathcal{\hat{K}}_{\Sigma}M$. Note a leaf is
preserved by the induced almost Hermitian structure on $\mu^{-1}\left(
0\right)  $ from the Hermitian complex structure $J$ on $M$, because a tangent
vector on it is a section to the subbundle $T_{\Sigma}+JT_{\Sigma}$. However
$\mathcal{\hat{K}}_{\Sigma}M$ may not be symplectic, because the symplectic
form $\omega_{Map}\;$on $\mu^{-1}\left(  0\right)  $ descends to a 2-form
$\omega^{\mathcal{K}}$ on $\mathcal{\hat{K}}_{\Sigma}M$ which may not be closed.

We will show that $\mathcal{\hat{K}}_{\Sigma}M$ admits three almost complex
structures $I^{\mathcal{K}}$, $K^{\mathcal{K}}$ and $J^{\mathcal{K}}$
satisfying the Hamilton relation%
\[
\left(  I^{\mathcal{K}}\right)  ^{2}=\left(  J^{\mathcal{K}}\right)
^{2}=\left(  J^{\mathcal{K}}\right)  ^{2}=I^{\mathcal{K}}J^{\mathcal{K}%
}K^{\mathcal{K}}=-id.
\]
Furthermore the associated K\"{a}hler forms $\omega_{I}^{\mathcal{K}}$ and
$\omega_{K}^{\mathcal{K}}$ are closed. In \cite{H2}, Hitchin showed that
existence of such structures on any finite dimensional Riemannian manifold
implies that its almost complex structure $J^{\mathcal{K}}$ is
\textit{integrable. }Namely we obtain a 1-fold $\mathbb{C}$-VCP on a Hermitian
complex manifold. If, in addition, $\omega_{J}^{\mathcal{K}}$ is closed then
$I^{\mathcal{K}}$, $J^{\mathcal{K}}$ and $K^{\mathcal{K}}$ are all integrable
complex structures and we have a hyperk\"{a}hler manifold.

\bigskip\ 

We begin with the following lemma which holds true both on $\mu^{-1}\left(
0\right)  $ and on $\mathcal{\hat{K}}_{\Sigma}M$.

\begin{lemma}
On $\mu^{-1}\left(  0\right)  \subset Map\left(  \Sigma,M\right)  $, the
endomorphisms $I$, $J$ and $K$ satisfy the following relations,
\begin{align*}
\;IJ  &  =-JI\text{ and }KJ=-JK,\\
\,I  &  =-KJ\;\text{and}\;K=IJ.
\end{align*}
\end{lemma}

Proof: The formula $IJ=-JI$ can be restated as follows: for any $f\in\mu
^{-1}\left(  0\right)  $ and for any tangent vectors of $Map\left(
\Sigma,M\right)  $ at $f$, $A,B$ $\in\Gamma(\Sigma,f^{\ast}\left(
T_{M}\right)  )$ we have
\[
g_{f}\left(  I_{f}J_{f}\left(  A\right)  ,B\right)  =g_{f}\left(  -J_{f}%
I_{f}\left(  A\right)  ,B\right)  .
\]
For simplicity, we ignore the subscript $f$. Since $J^{2}=-id$, we have
\begin{align*}
&  g\left(  IJ\left(  A\right)  ,B\right)  -g\left(  -JI\left(  A\right)
,B\right) \\
&  =g\left(  IJ\left(  A\right)  ,B\right)  -g\left(  I\left(  A\right)
,JB\right) \\
&  =\omega_{I}\left(  JA,B\right)  -\omega_{I}\left(  A,JB\right) \\
&  =\int_{\Sigma}\iota_{JA\wedge B}\operatorname{Re}\,\Omega_{M}-\int_{\Sigma
}\iota_{A\wedge JB}\operatorname{Re}\,\Omega_{M}\\
&  =\int_{\Sigma}\iota_{(JA\wedge B-A\wedge JB)}\frac{1}{2}(\,\Omega_{M}%
+\bar{\Omega}_{M})=0\text{.}%
\end{align*}
The last equality follows from the fact that $\Omega_{M}\;$and $\bar{\Omega
}_{M}\;$vanish when each is contracted by an element in $\wedge^{1,1}T_{M}$.
By replacing $\operatorname{Re}\,\Omega_{M}$ with $\operatorname{Im}%
\,\Omega_{M}$ in the above calculations, we also have $KJ=-JK$.

To prove the others formulas, we consider%

\begin{align*}
&  g(IA,B)-g(-KJA,B)\\
&  =\omega_{I}(A,B)+\omega_{K}(JA,B)\\
&  =\int_{\Sigma}\iota_{A\wedge B}\operatorname{Re}\,\Omega_{M}+\int_{\Sigma
}\iota_{JA\wedge B}\left(  -\operatorname{Im}\,\Omega_{M}\right) \\
&  =\int_{\Sigma}\iota_{A\wedge B}\frac{1}{2}\left(  \Omega_{M}+\bar{\Omega
}_{M}\right)  -\int_{\Sigma}\iota_{JA\wedge B}\frac{-\sqrt{-1}}{2}\left(
\Omega_{M}-\bar{\Omega}_{M}\right) \\
&  =\frac{1}{2}\int_{\Sigma}\iota_{(A+\sqrt{-1}JA)\wedge B}\,\Omega_{M}%
+\frac{1}{2}\int_{\Sigma}\iota_{(A-\sqrt{-1}JA)\wedge B}\,\bar{\Omega}_{M}=0
\end{align*}
since $i_{(A+\sqrt{-1}JA)}\,\Omega_{M}\;=0\;$and taking complex conjugation.
This implies that $g(IA,B)-g(-KJA,B)=0$ and therefore $I=-KJ$. Finally
$J^{2}=-Id$ and $I=-KJ$ imply that $K=IJ.$ Hence the results. $\blacksquare$

\begin{theorem}
Suppose that $M$ is a Calabi-Yau n-fold $M$. For any $n-2$ dimensional closed
manifold $\Sigma$ with $b_{1}\left(  \Sigma\right)  =0$, the isotropic knot
space $\mathcal{\hat{K}}_{\Sigma}M$ is an infinite dimensional integrable
complex manifold with a natural 1-fold $\mathbb{C}$-VCP structure, in
particular a natural holomorphic symplectic structure.
\end{theorem}

\begin{proof}
From the construction of $\mathcal{\hat{K}}_{\Sigma}M$, it has a Hermitian
K\"{a}hler form $\omega^{\mathcal{K}}$ induced from that of $M$ and a closed
holomorphic symplectic form $\Omega^{\mathcal{K}}$ given by the transgression
of the closed $(n-1)$-fold $\mathbb{C}$-\textit{VCP} form $\Omega_{M}$ on $M$.
As we have seen above, the induced holomorphic symplectic form is closed but
the induced Hermitian K\"{a}hler form may not be closed. If their
corresponding almost complex structures satisfy the Hamilton relation then
this implies that $J^{\mathcal{K}}$ is integrable \cite{H2}. Namely
$\mathcal{\hat{K}}_{\Sigma}M$ is a Hermitian integrable complex manifold with
a 1-fold $\mathbb{C}$-VCP structure. In order to verify the Hamilton
relation,
\[
\left(  I^{\mathcal{K}}\right)  ^{2}=\left(  J^{\mathcal{K}}\right)
^{2}=\left(  J^{\mathcal{K}}\right)  ^{2}=I^{\mathcal{K}}J^{\mathcal{K}%
}K^{\mathcal{K}}=-id,
\]
we only need to show that $\left(  I^{\mathcal{K}}\right)  ^{2}=-Id$ and
$\left(  K^{\mathcal{K}}\right)  ^{2}=$ $-Id$ because of the previous lemma.
Namely, $I^{\mathcal{K}}$ and $K^{\mathcal{K}}$ are almost complex structures
on $\mathcal{\hat{K}}_{\Sigma}M$.

We consider a fixed $\left[  f\right]  $ in $\mathcal{\hat{K}}_{\Sigma}M$ and
by localization method as in the proof of lemma \ref{Localize Lemma}, we can
reduce the identities to the tangent space of a point $x$ in $\Sigma$. The
transgression
\[
\int_{\Sigma}ev^{\ast}\Omega_{M}=\omega_{I,\left[  f\right]  }^{\mathcal{K}%
}-\sqrt{-1}\omega_{K,\left[  f\right]  }^{\mathcal{K}}%
\]
is descended to
\[
\iota_{E_{1}\wedge...\wedge E_{n-2}}\,\Omega_{M}=\omega_{I,\left[  f\right]
,x}^{\mathcal{K}}-\sqrt{-1}\omega_{K,\left[  f\right]  ,x}^{\mathcal{K}}%
\]
where $E_{1,}...,E_{n-2}$ is an orthonormal oriented basis $\left(  T_{\Sigma
}\right)  _{x}$. Since $f$ is isotropic, the complexified vectors $E_{i}%
-\sqrt{-1}JE_{i}$ of $E_{i}$ can be defined over $\left(  T_{\Sigma
}+JT_{\Sigma}\right)  _{x}$, therefore the above equality is equivalent to
\[
\iota_{\left(  E_{1}-\sqrt{-1}JE_{1}\right)  /2\wedge...\wedge\left(
E_{n-2}-\sqrt{-1}JE_{n-2}\right)  /2}\,\Omega_{M}=\omega_{I,\left[  f\right]
,x}^{\mathcal{K}}-\sqrt{-1}\omega_{K,\left[  f\right]  ,x}^{\mathcal{K}},
\]
i.e. a $1$-fold $\mathbb{C}$-\textit{VCP} on $f^{\ast}\left(  TM\right)
/\left(  T_{\Sigma}+JT_{\Sigma}\right)  \mid_{x}$ which is $T_{\left[
f\right]  ,x}\mathcal{\hat{K}}_{\Sigma}M$. Since this $1$-fold $\mathbb{C}%
$-\textit{VCP }gives a hyperk\"{a}hler structure on $T_{\left[  f\right]
,x}\mathcal{\hat{K}}_{\Sigma}M$, $I_{\left[  f\right]  }^{\mathcal{K}}$ and
$K_{\left[  f\right]  }^{\mathcal{K}}$ satisfy the Hamilton relation at $x$ in
$\Sigma$. That is
\[
\left(  I_{\left[  f\right]  ,x}^{\mathcal{K}}\right)  ^{2}=-Id\text{ and
}\left(  K_{\left[  f\right]  ,x}^{\mathcal{K}}\right)  ^{2}=-Id.
\]
Therefore we have $\left(  I^{\mathcal{K}}\right)  ^{2}=-Id$ and $\left(
K^{\mathcal{K}}\right)  ^{2}=$ $-Id$. Hence the result.
\end{proof}

\begin{remark}
\label{CY3 remark}In String theory we need to compactify a ten dimensional
spacetime on a Calabi-Yau threefold. When $M$ is a Calabi-Yau threefold, then
$\Sigma$ is an one dimensional circle and therefore $b_{1}\left(
\Sigma\right)  $ is nonzero. In general, as discussed in \cite{D} and
\cite{H1}, the symplectic quotient construction for $Map\left(  \Sigma
,M\right)  //Diff\left(  \Sigma,\sigma\right)  $ can be modified to obtain a
symplectic structure on a rank $b_{1}\left(  \Sigma\right)  $ torus bundle
over the space of isotropic submanifolds in $M$. Roughly speaking this torus
bundle is the space of isotropic submanifolds coupled with flat rank one line
bundles (or gerbes) in $M$. In the Calabi-Yau threefold case, every circle
$\Sigma$ in $M$ is automatically isotropic. Therefore $\mathcal{\hat{K}%
}_{\Sigma}M$ is the space of loops (or \emph{string}) coupled with flat line
bundles in $M$, up to deformations of strings along their complexified tangent
directions. We wonder whether this infinite dimensional holomorphic symplectic
manifold $\mathcal{\hat{K}}_{\Sigma}M$ has any natural physical interpretations.
\end{remark}

\subsection{\label{Sec CLag Isot Knot}Complex Lagrangians in Isotropic Knot Spaces}

In this subsection we relate the geometry of $\mathbb{C}$-VCP of a Calabi-Yau
manifold $M$ to the holomorphic symplectic geometry of its isotropic knot
space $\mathcal{\hat{K}}_{\Sigma}M$. For example both N-branes (i.e. complex
hypersurfaces) and D-branes (i.e. special Lagrangian submanifold with phase
$-\pi/2$) in $M$ correspond to complex Lagrangian submanifold in
$\mathcal{\hat{K}}_{\Sigma}M$, but for different almost complex structures in
the twistor family. First we discuss the correspondence for instantons.

In the following proposition, we use $e^{i\theta}\Omega_{M}$\ instead of
$\Omega_{M}$ to get a $1$-fold $\mathbb{C}$-\textit{VCP} on $\mathcal{\hat{K}%
}_{\Sigma}M$, and also $\operatorname{Re}\left(  e^{i\theta}\Omega_{M}\right)
$ gives another symplectic form $\mathbb{\omega}_{I,\theta}^{\mathcal{K}}%
\;$and corresponding Hermitian almost complex structure $J_{\theta
}^{\mathcal{K}}\;$on $\mathcal{\hat{K}}_{\Sigma}M$ defined as $J_{\theta
}^{\mathcal{K}}= $ $\cos\theta I^{\mathcal{K}}+\cos\theta K^{\mathcal{K}}$.

\begin{proposition}
Suppose that $M$ is a Calabi-Yau $n$-fold. Let $D$ be a normal disk$\;$in
$\mu^{-1}\left(  0\right)  $ with an n-dimensional submanifold $A$ defined
as$\;A=D\tilde{\times}\Sigma$ , and assume $A\rightarrow D$ is a Riemannian
submersion. We denote the reduction$\;$of $D\;$in$\;\mathcal{\hat{K}}_{\Sigma
}M$ as $\hat{D}$.

Then $\hat{D}$ is an instanton$\;$i.e. a $J_{\theta}^{\mathcal{K}}%
$-holomorphic curve in $\mathcal{\hat{K}}_{\Sigma}M$ if and only if $A\;$is an
instanton i.e. a special Lagrangian with phase $\theta$ in $M$.
\end{proposition}

\begin{proof}
In the proposition \ref{ProDD}, VCP form $\phi$ plays a role of calibrating
form rather than that of a VCP form. So by replacing $\phi$ by
$\operatorname{Re}\left(  e^{i\theta}\Omega_{M}\right)  $ in proposition
\ref{ProDD}, readers can see this theorem can be proved in essentially the
same manner. But since proposition \ref{ProDD} is given for the parametrized
knot space and this theorem is given by a reduction, we need to check $\hat
{D}$ is a disk in $\mathcal{\hat{K}}_{\Sigma}M$. Let $f$ be the center of disk
$D.$ Since $D$ is a normal disk, a tangent vector $v$ at $f$ along
$D\subset\mu^{-1}\left(  0\right)  $ is in $\Gamma\left(  \Sigma,N_{\Sigma
/M}\right)  $. And since $f$ is isotropic, $Jv$ is also in $\Gamma\left(
\Sigma,N_{\Sigma/M}\right)  $ equivalently $v$ is in $\Gamma\left(
\Sigma,N_{J\Sigma/M}\right)  $. So $v$ $\in\Gamma\left(  \Sigma,f^{\ast
}\left(  TM\right)  /\left(  T_{\Sigma}+JT_{\Sigma}\right)  \right)  $. This
implies $T_{f}$ $D$ can be identified with $T_{\left[  f\right]  }\hat{D}$,
i.e. $\hat{D}$ is a disk in $\mathcal{\hat{K}}_{\Sigma}M$.
\end{proof}

\bigskip

Next we are going to relate N- and D-branes in Calabi-Yau manifolds $M$ to
complex Lagrangian submanifolds in $\mathcal{\hat{K}}_{\Sigma}M$.

\begin{definition}
For any submanifold $C$ in a Calabi-Yau manifold $M$, we define a subspace
$\mathcal{\hat{K}}_{\Sigma}C$ in $\mathcal{\hat{K}}_{\Sigma}M=\mu^{-1}\left(
0\right)  /\left\langle \mathcal{D}\right\rangle $ as follow: The equivalent
relation $\left\langle \mathcal{D}\right\rangle $ on $\mu^{-1}\left(
0\right)  $ restrict to one on $Map\left(  \Sigma,C\right)  \cap\mu
^{-1}\left(  0\right)  $ and we define
\[
\mathcal{\hat{K}}_{\Sigma}C=\left\{  Map\left(  \Sigma,C\right)  \cap\mu
^{-1}\left(  0\right)  \right\}  /\left\langle \mathcal{D}\right\rangle
\text{.}%
\]
\end{definition}

In the following theorem, we see that $C\;$being a $N$-brane in $M$
corresponds to $\mathcal{\hat{K}}_{\Sigma}C$ being a $J^{\mathcal{K}}$-complex
Lagrangian in $\mathcal{\hat{K}}_{\Sigma}M$ which means that $\mathcal{\hat
{K}}_{\Sigma}C$ is both maximally self $\omega_{I}^{\mathcal{K}}%
$-perpendicular and maximally self $\omega_{K}^{\mathcal{K}}$-perpendicular in
$\mathcal{\hat{K}}_{\Sigma}M$.

\begin{theorem}
Let $C$ be a connected analytic submanifold in a Calabi-Yau manifold $M$. Then
the following two statements are equivalent:

(i) $C$ is an $N$-brane (i.e. a complex hypersurface) in $M$;

(ii) $\mathcal{\hat{K}}_{\Sigma}C$ is a $J^{\mathcal{K}}$-complex Lagrangian
submanifold in $\mathcal{\hat{K}}_{\Sigma}M$.
\end{theorem}

\textbf{Proof: }For the \textit{if }part, we assume that $\mathcal{\hat{K}%
}_{\Sigma}C\;$is a $J^{\mathcal{K}}$-complex Lagrangian$\;$in $\mathcal{\hat
{K}}_{\Sigma}M$. First we need to show that the dimension of $C$ is at least
$2n-2$ where $n$ is the complex dimension of $M$.

Before proving this claim for the general case, let us discuss a simplified
linear setting where some of the key arguments become more transparent.
Suppose that $M$ is a \textit{linear} Calabi-Yau manifold, that is
$M\cong\mathbb{C}^{n}$ with the standard K\"{a}hler structure and $\Omega
_{M}=dz^{1}\wedge dz^{2}\wedge\cdots\wedge dz^{n}$. Let $\Sigma$ be a $\left(
n-2\right)  $-dimensional isotropic linear subspace in $M$ lying inside
another linear subspace $C$ in $M$. For simplicity we assume that $\Sigma$ is
the linear span of $x^{1},x^{2},...,x^{n-2}$. Of course $M/\Sigma\oplus
J\Sigma\cong\mathbb{C}^{2}$ is a linear holomorphic symplectic manifold with
\[
\int_{\Sigma}\Omega_{M}=dz^{n-1}\wedge dz^{n},
\]
which is the standard holomorphic symplectic form on $\mathbb{C}^{2}$. Suppose
that
\[
C/\left(  \Sigma\oplus J\Sigma\right)  \cap C\subset M/\Sigma\oplus J\Sigma
\]
is a complex Lagrangian subspace. Then there is a vector in $C$ perpendicular
to both $\Sigma$ and $J\Sigma$, say $x^{n-1}$. If we denote the linear span of
$x^{2},...,x^{n-2},x^{n-1}$ as $\Sigma^{\prime}$, then $\Sigma^{\prime}$ is
another $\left(  n-2\right)  $-dimensional isotropic linear subspace in $M$
lying inside $C$. Furthermore $x^{1}$ is a normal vector in $C$ perpendicular
to $\Sigma^{\prime}\oplus J\Sigma^{\prime}$. This implies that $y^{1}=Jx^{1}$
also lie in $C$. This is because $C/\left(  \Sigma^{\prime}\oplus
J\Sigma^{\prime}\right)  \cap C\subset M/\Sigma^{\prime}\oplus J\Sigma
^{\prime}$ being a complex Lagrangian subspace implies that it is invariant
under the complex structure on $M/\Sigma^{\prime}\oplus J\Sigma^{\prime}$
induced by $J$ on $M$. By the same reasoning, $y^{j}$ also lie in $C$ for
$j=1,2,...,n-2$. That is $C$ contains the linear span of $\left\{  x^{j}%
,y^{j}\right\}  _{j=1}^{n-2}$. On the other hand it also contain $x^{n-1}$ and
$y^{n-1}$ and therefore $\dim C\geq2n-2$.

We come back to the general situation where $\mathcal{\hat{K}}_{\Sigma}C\;$is
a $J^{\mathcal{K}}$-complex Lagrangian$\;$in $\mathcal{\hat{K}}_{\Sigma}M $.
One difficulty is to rotate the isotropic submanifold $\Sigma$ to
$\Sigma^{\prime}$ inside $M$.

We observe that the tangent space to $\mathcal{\hat{K}}_{\Sigma}C$ at any
point $\left[  f\right]  $ is given by
\[
T_{\left[  f\right]  }\left(  \mathcal{\hat{K}}_{\Sigma}C\right)
=\Gamma\left(  \Sigma,\frac{f^{\ast}T_{C}}{T_{\Sigma}\oplus JT_{\Sigma}\cap
f^{\ast}T_{C}}\right)  \text{.}%
\]
Since any $J^{\mathcal{K}}$-complex Lagrangian submanifold is indeed an
integrable complex submanifold, $T_{\left[  f\right]  }\left(  \mathcal{\hat
{K}}_{\Sigma}C\right)  $ is preserved by $J^{\mathcal{K}}$. This implies that
the complex structure on $T_{M}$ induces a complex structure on the quotient
bundle $f^{\ast}T_{C}/\left[  T_{\Sigma}\oplus JT_{\Sigma}\cap f^{\ast}%
T_{C}\right]  $. For any $\nu\in T_{\left[  f\right]  }\left(  \mathcal{\hat
{K}}_{\Sigma}C\right)  $ we can regard it as a section of $f^{\ast}T_{C}$ over
$\Sigma$, perpendicular to $T_{\Sigma}$ and $JT_{\Sigma}\cap f^{\ast}T_{C}$.
In particular $Jv$ is also such a section.

We need the following lemma which will be proven later.

\begin{lemma}
In the above situation, we have $JT_{\Sigma}\subset f^{\ast}T_{C}$.
\end{lemma}

This lemma implies that the restriction of $T_{C}$ on $\Sigma$ contains the
linear span of $T_{\Sigma}$, $JT_{\Sigma}$, $\nu$ and $J\nu\,$, which has rank
$2n-2$. (More precisely we consider the restriction to the complement of the
zero set of $\nu$ in $\Sigma$.) Therefore $\dim C\geq2n-2$.

On the other hand, if $\dim C$ $>2n-2$, then by using localizing method as in
lemma \ref{Localize Lemma}, $\mathcal{\hat{K}}_{\Sigma}C$ would be too large
to be a complex Lagrangian submanifold in $\mathcal{\hat{K}}_{\Sigma}M$. That
is $\dim C=2n-2$.

In particular, this implies that $f^{\ast}T_{C}$ is isomorphic to the linear
span of $T_{\Sigma}$, $JT_{\Sigma}$, $\nu$ and $J\nu$ outside the zero set of
$\nu$. Therefore $T_{C}$ is preserved by the complex structure of $M$ along
$\Sigma$. By varying the isotropic submanifold $\Sigma$ in $C$, we can cover
an open neighborhood of $\Sigma$ in $C$ (for example using the gluing
arguments as in the proof of the above lemma). This implies that an open
neighborhood of $\Sigma$ in $C$ is a complex submanifold in $M$. By the
analyticity of $C$, the submanifold $C$ is a complex hypersurface in $M$.

For the \textit{only if }part, we suppose $C$ is complex hypersurface in $M$,
then it is clear that $\omega_{I}^{\mathcal{K}}$ and $\omega_{K}^{\mathcal{K}%
}$ vanish along $\mathcal{\hat{K}}_{\Sigma}C$ since $\Omega_{M}$ vanishes
along $C$. Using similar arguments as in the proof of proposition \ref{ProD},
it is not difficult to verify that $\mathcal{\hat{K}}_{\Sigma}C$ is maximally
self $\omega_{I}^{\mathcal{K}}$-perpendicular and maximally self $\omega
_{K}^{\mathcal{K}}$-perpendicular in $\mathcal{\hat{K}}_{\Sigma}M$. That is
$\mathcal{\hat{K}}_{\Sigma}C$ is a $J^{\mathcal{K}}$-complex Lagrangian
submanifold in $\mathcal{\hat{K}}_{\Sigma}M$. Hence the theorem. $\blacksquare$

\bigskip

We suspect that the analyticity assumption on $C$ is unnecessary. All we need
is to deform the isotropic submanifold $\Sigma$ of $M$ inside $C$ to cover
every point in $C$.

\bigskip

\textbf{Proof of lemma}: For any tangent vector $u$ of $\Sigma$ at a point $p
$, i.e. $u\in T_{p}\Sigma\subset T_{p}C$, we need to show that $Ju\in T_{p}C$.
We can assume that $\nu\left(  p\right)  $ has unit length. First we choose
local holomorphic coordinate $z^{1},z^{2},...,z^{n}$ near $p$ satisfying the
following properties at $p$:
\begin{align*}
z^{i}\left(  p\right)   &  =0\text{ for all }i=1,...,n,\\
\Omega_{M}\left(  p\right)   &  =dz^{1}\wedge dz^{2}\wedge\cdots\wedge
dz^{n},\\
&  T_{p}\Sigma\text{ is spanned by }\frac{\partial}{\partial x^{1}}%
,\frac{\partial}{\partial x^{2}},...,\frac{\partial}{\partial x^{n-2}},\\
\nu\left(  p\right)   &  =\frac{\partial}{\partial x^{n-1}}\text{ and }%
u=\frac{\partial}{\partial x^{1}}\text{.}%
\end{align*}
where $z^{j}=x^{j}+iy^{j}$ for each $j$. We could use $x^{1},...,x^{n-2}$ to
parametrize $\Sigma$ near $p$.

Recall that
\[
\nu\in T_{\left[  f\right]  }\left(  \mathcal{\hat{K}}_{\Sigma}C\right)
=\Gamma\left(  \Sigma,\frac{f^{\ast}T_{C}}{T_{\Sigma}\oplus JT_{\Sigma}\cap
f^{\ast}T_{C}}\right)  \text{.}%
\]
For any smooth function $\alpha\left(  x^{1},...,x^{n-2}\right)  $ on $\Sigma
$, $\alpha\nu$ is again a tangent vector to $\mathcal{\hat{K}}_{\Sigma}C$ at
$\left[  f\right]  $. We are going to construct a particular $\alpha$
supported on a small neighborhood of $p$ in $\Sigma$ satisfying
\begin{align*}
\frac{\partial\alpha}{\partial x^{1}}\left(  0\right)   &  =\infty,\\
\frac{\partial\alpha}{\partial x^{j}}\left(  0\right)   &  =0\text{, for
}j=2,...,n-2\text{.}%
\end{align*}
For simplicity we assume that this small neighborhood contain the unit ball in
$\Sigma$. To construct $\alpha$, we fix a smooth even cutoff function
$\eta:\mathbb{R}\rightarrow\left[  0,1\right]  $ with $supp\left(
\eta\right)  \in\left(  -1,1\right)  $, $\eta\left(  0\right)  =1$ and
$\eta^{\prime}\left(  0\right)  =0$. We define $\alpha:\Sigma\rightarrow
\mathbb{R}$ as follow:
\[
\alpha\left(  x^{1},x^{2},...,x^{n-2}\right)  =\left(  x^{1}\right)
^{1/3}\cdot\eta\left(  x^{1}\right)  \cdot\eta\left(
{\textstyle\sum_{j=2}^{n-2}}
\left(  x^{j}\right)  ^{2}\right)  \text{.}%
\]

For any small real number $\varepsilon$, we write
\[
\nu_{\varepsilon}=\varepsilon\alpha\nu\in T_{\left[  f\right]  }\left(
\mathcal{\hat{K}}_{\Sigma}C\right)
\]
and we denote the corresponding family of isotropic submanifolds of $M$ in $C
$ by $f_{\varepsilon}:\Sigma_{\varepsilon}\rightarrow M$. From $\alpha\left(
0\right)  =0$, there is a family of points $p_{\varepsilon}\in\Sigma
_{\varepsilon}$ with the property that
\[
\left|  p_{\varepsilon}-p\right|  =O\left(  \varepsilon^{2}\right)  \text{.}%
\]

Using the property $\partial\alpha/\partial x^{1}\left(  0\right)  =\infty$,
we can find a family of normal vectors at $p_{\varepsilon},$
\[
u_{\varepsilon}\in N_{\Sigma_{\varepsilon}/C},
\]
with the property
\[
\left|  u_{\varepsilon}-\frac{\partial}{\partial x^{1}}\right|  =O\left(
\varepsilon^{2}\right)  \text{.}%
\]

Since $T_{\left[  f_{\varepsilon}\right]  }\left(  \mathcal{\hat{K}}_{\Sigma
}C\right)  $ is preserved by $J^{\mathcal{K}}$, $f_{\varepsilon}^{\ast}%
T_{C}/\left[  T_{\Sigma_{\varepsilon}}\oplus JT_{\Sigma_{\varepsilon}}\cap
f_{\varepsilon}^{\ast}T_{C}\right]  $ is preserved by $J$, we have
$Ju_{\varepsilon}\in T_{p_{\varepsilon}}C$. By letting $\varepsilon$ goes to
zero, we have $Ju\in T_{p}C$. Hence the lemma. $\blacksquare$

\bigskip

In the following theorem, we see that relationship between a $D$-brane $L$\ in
$M\;$with respect to $e^{i\theta}\Omega_{M}$ and a $J_{\theta+\frac{\pi}{2}%
}^{\mathcal{K}}$ complex Lagrangian $\mathcal{\hat{K}}_{\Sigma}L$ in
$\mathcal{\hat{K}}_{\Sigma}M$, i.e. it is maximally self $\omega_{I,\theta
}^{\mathcal{K}}$-perpendicular and maximally self $\omega^{\mathcal{K}}$-perpendicular$.$

\begin{theorem}
Let $L\;$be a connected analytic submanifold of a Calabi-Yau manifold $M$,
then the following two statements are equivalent:

(i) $L$ is a \textit{D}-brane with phase $\theta,$(i.e. a special Lagrangian
with phase $\theta-\frac{\pi}{2}$);

(ii)$\;\mathcal{\hat{K}}_{\Sigma}L$ is a $J_{\theta+\frac{\pi}{2}%
}^{\mathcal{K}}$-complex Lagrangian submanifold in $\mathcal{\hat{K}}_{\Sigma
}M$.
\end{theorem}

\textbf{Proof}: By replacing the holomorphic volume form on $M$ from
$\Omega_{M}$ to $e^{i\theta}\Omega_{M}$ if necessary, we can assume that
$\theta$ is zero.

For the \textit{only if }part, we assume that $L$ is a special Lagrangian
submanifold in $M$ with phase $-\pi/2$. Because the K\"{a}hler form
$\mathbb{\omega}$ and $\operatorname{Re}\Omega_{M}$ of $M$ vanish along $L$,
it is clear that $\mathbb{\omega}^{\mathcal{K}}$ and $\mathbb{\omega}%
_{I,0}^{\mathcal{K}}$ vanishes along $\mathcal{\hat{K}}_{\Sigma}L$. Notice
that $L$ being a Lagrangian submanifold in $M$ implies that any submanifold in
$L$ is automatically isotropic in $M$, and moreover the equivalent relation
$\left\langle \mathcal{D}\right\rangle $ on $\mu^{-1}\left(  0\right)  $
restricted to $Map\left(  \Sigma,L\right)  $ is trivial locally. That is
$\mathcal{\hat{K}}_{\Sigma}L$ is the same as $\mathcal{K}_{\Sigma}L=Map\left(
\Sigma,L\right)  /Diff\left(  \Sigma\right)  $ at least locally.

We claim that $\mathcal{\hat{K}}_{\Sigma}L$ is maximally self $\omega
^{\mathcal{K}}$-perpendicular in $\mathcal{\hat{K}}_{\Sigma}M$. Otherwise
there is normal vector field $\nu\in\Gamma\left(  \Sigma,N_{\Sigma/M}\right)
$ not lying in $\Gamma\left(  \Sigma,N_{\Sigma/L}\right)  $ such that
$\omega\left(  \nu,u\right)  =0$ for any $u\in\Gamma\left(  \Sigma,f^{\ast
}T_{L}\right)  $. Suppose that $\nu\left(  p\right)  \notin N_{\Sigma/L,p}$,
then this implies that $\omega$ vanishes on the linear span of $T_{p}L$ and
$\nu\left(  p\right)  $ inside $T_{p}M$. This is impossible because the
dimension of the linear span is bigger than $n$.

Similarly $\mathcal{\hat{K}}_{\Sigma}L$ is maximally self $\omega
_{I,0}^{\mathcal{K}}$-perpendicular in $\mathcal{\hat{K}}_{\Sigma}M$.
Otherwise $\operatorname{Re}\Omega_{M}$ vanishes on the linear span of
$T_{p}L$ and $\nu\left(  p\right)  $ inside $T_{p}M$. This is again impossible
because $\Omega_{M,p}$ is a complex volume form on $T_{p}M\cong\mathbb{C}^{n}$
and therefore can not vanish on any co-isotropic subspaces other than
Lagrangians. Hence $\mathcal{\hat{K}}_{\Sigma}L\;$ is a $J_{\pi/2}%
^{\mathcal{K}}$-complex Lagrangian in $\mathcal{\hat{K}}_{\Sigma}M$.

\bigskip

For the \textit{if} part, we assume that $\mathcal{\hat{K}}_{\Sigma}L\;$ is a
$K^{\mathcal{K}}$-complex Lagrangian submanifold in $\mathcal{\hat{K}}%
_{\Sigma}M$. For any $\left[  f\right]  \in\mathcal{\hat{K}}_{\Sigma}L$, the
tangent space of $\mathcal{\hat{K}}_{\Sigma}L$ (resp. $\mathcal{\hat{K}%
}_{\Sigma}M$) at $\left[  f\right]  $ is the section of the bundle $f^{\ast
}T_{L}/T_{\Sigma}\oplus JT_{\Sigma}\cap f^{\ast}T_{L}$ (resp. $f^{\ast}\left(
T_{M}\right)  /\left(  T_{\Sigma}+JT_{\Sigma}\right)  $) over $\Sigma$. Let
$v\in T_{\left[  f\right]  }\left(  \mathcal{\hat{K}}_{\Sigma}L\right)  $, we
can regard $v$ as a section of $f^{\ast}T_{L}$ over $\Sigma$, orthogonal to
$T_{\Sigma}\oplus JT_{\Sigma}\cap f^{\ast}T_{L}$.

Since $f^{\ast}\left(  T_{M}\right)  /\left(  T_{\Sigma}+JT_{\Sigma}\right)  $
is a rank four bundle and $\mathcal{\hat{K}}_{\Sigma}L$ is a Lagrangian in
$\mathcal{\hat{K}}_{\Sigma}M$, this implies that $f^{\ast}T_{L}/T_{\Sigma
}\oplus JT_{\Sigma}\cap f^{\ast}T_{L}$ must be a rank two bundle over $\Sigma
$. Therefore
\[
\dim L\geq n,
\]
with the equality sign holds if and only if $JT_{\Sigma}\cap f^{\ast}T_{L}$ is
trivial. Suppose that $\nu$ is a tangent vector of $\mathcal{\hat{K}}_{\Sigma
}L$ at $f$

Assume that $JT_{\Sigma}\cap f^{\ast}T_{L}$ is not trivial, we can find a
tangent vector $u$ to $f\left(  \Sigma\right)  $ at a point $p$ such that
$Ju\in T_{p}L$. For simplicity we assume that $\nu$ has unit length at the
point $p$.

As in the proof of the previous theorem, we can choose local holomorphic
coordinates $z^{j}$'s of $M$ around $p$ such that%

\begin{align*}
z^{i}\left(  p\right)   &  =0\text{ for all }i=1,...,n,\\
\Omega_{M}\left(  p\right)   &  =dz^{1}\wedge dz^{2}\wedge\cdots\wedge
dz^{n},\\
&  T_{p}\Sigma\text{ is spanned by }\frac{\partial}{\partial x^{1}}%
,\frac{\partial}{\partial x^{2}},...,\frac{\partial}{\partial x^{n-2}}%
\text{,}\\
v\left(  p\right)   &  =\frac{\partial}{\partial x^{n-1}}\text{, }%
u=\frac{\partial}{\partial x^{1}}\text{and }Ju=\frac{\partial}{\partial y^{1}%
}\text{.}%
\end{align*}

We choose a function $\alpha\left(  x\right)  $ as in the proof of the
previous theorem, write
\[
\nu_{\varepsilon}=\varepsilon\alpha\nu\in T_{\left[  f\right]  }\left(
\mathcal{\hat{K}}_{\Sigma}L\right)
\]
for any small real number $\varepsilon$, and we denote the corresponding
family of isotropic submanifolds of $M$ in $L$ by $f_{\varepsilon}%
:\Sigma_{\varepsilon}\rightarrow M$ as before. From $\alpha\left(  0\right)
=0$, there is a family of points $p_{\varepsilon}\in\Sigma_{\varepsilon}$ with
the property that
\[
\left|  p_{\varepsilon}-p\right|  =O\left(  \varepsilon^{2}\right)  \text{.}%
\]

Using the property $\partial\alpha/\partial x^{1}\left(  0\right)  =\infty$,
we can find two family of normal vectors at $p_{\varepsilon},$
\[
u_{\varepsilon},w_{\varepsilon}\in N_{\Sigma_{\varepsilon}/L},
\]
with the property
\[
\left|  u_{\varepsilon}-\frac{\partial}{\partial x^{1}}\right|  =O\left(
\varepsilon^{2}\right)  \text{ and }\left|  w_{\varepsilon}-\frac{\partial
}{\partial y^{1}}\right|  =O\left(  \varepsilon^{2}\right)  \text{.}%
\]

However, using localization arguments as before, $\mathcal{\hat{K}}_{\Sigma}L$
being an $\omega^{\mathcal{K}}$-Lagrangian submanifold in $\mathcal{\hat{K}%
}_{\Sigma}M$ implies that $\omega\left(  u_{\varepsilon},w_{\varepsilon
}\right)  =0$. By letting $\varepsilon$ goes to zero, we have $0=\omega\left(
\partial/\partial x^{1},\partial/\partial y^{1}\right)  =-1$, a contradiction.
Hence $JT_{\Sigma}\cap f^{\ast}T_{L}$ is trivial and $\dim L=n$.\newline 

Since $\Gamma\left(  \Sigma,f^{\ast}T_{L}/T_{\Sigma}\oplus JT_{\Sigma}\cap
f^{\ast}T_{L}\right)  $ is preserved by $K^{\mathcal{K}}$, we have
$K^{\mathcal{K}}\nu\left(  p\right)  \in f^{\ast}T_{L}$. Note that
$K^{\mathcal{K}}\nu\left(  p\right)  $ is the tangent vector of $M$ which is
the metric dual of the one form,
\begin{align*}
\nu\lrcorner\omega_{K}^{\mathcal{K}}\left(  p\right)   &  =\left(
\frac{\partial}{\partial x^{1}}\wedge\frac{\partial}{\partial x^{2}}%
\wedge\cdots\wedge\frac{\partial}{\partial x^{n-2}}\wedge\nu\right)
\lrcorner\operatorname{Im}\Omega_{M}\left(  p\right) \\
&  =\pm dy^{n}\text{,}%
\end{align*}
since $\nu\left(  p\right)  =\partial/\partial x^{n-1}$. This implies that
$T_{p}L$ is spanned by $\partial/\partial x^{j}$'s for $j=1,2,...,n-1$ and
$\partial/\partial y^{n}$, that is a special Lagrangian subspace of phase
$-\pi/2$.

As in the proof of the previous theorem, by deforming the isotropic
submanifold $\Sigma$ in $L$ and using the fact that $L$ is an analytic
submanifold of $M$, we conclude that $L$ is a special Lagrangian submanifold
in $M$ with zero phase.

Hence the theorem. $\blacksquare$

\section{Concluding Remarks}

In this paper we study both real and complex vector cross products (VCP).
Instantons in either settings are calibrated submanifolds. This gives a
unified way to explain the calibrating property of many such examples, as
studied by Harvey and Lawson in \cite{H-L}. It is desirable to further study
the calibration geometry from this point of view.

Manifolds with real VCP include symplectic/K\"{a}hler and $G_{2}$-manifolds.
In section \ref{Symplectic on Loop}, we relate the geometry of VCP on the
manifolds $M$ to the symplectic geometry of their knot spaces $\mathcal{K}%
_{\Sigma}M$. Motivated from this relationship, it is natural to study the
intersection theory of branes and count the number of instantons bounding
them, similar to the Floer's homology theory of Lagrangian intersections. For
example, in the case of $G_{2}$-manifolds, counting associative submanifolds
bounding nearby coassociative submanifolds is closely related to the
Seiberg-Witten invariants of the four dimensional coassociative submanifolds
\cite{Leung XWWang}. Results along this line should be useful in understanding
the M-theory in Physics.

Another interesting problem for VCP is the study of the submanifold geometry
of branes. For example one would like to have a unified approach to the mean
curvature flow for both hypersurfaces and Lagrangian submanifolds.

\bigskip

For manifolds with $\mathbb{C}$-VCP $M$, we classify them in theorem
\ref{Classification of CVCP}, namely they must be either holomorphic volume
forms in Calabi-Yau manifolds or holomorphic symplectic forms in
hyperk\"{a}hler manifolds. We study the geometry of instantons, Dirichlet
branes and Neumann branes in $M$. In section \ref{SecConstrModifLoop}, we
construct an isotropic knot space $\mathcal{\hat{K}}_{\Sigma}M$ for any
Calabi-Yau manifold $M$ and show that it admits a natural holomorphic
symplectic structure. We also relate the Calabi-Yau geometry of $M$ to the
holomorphic symplectic geometry of $\mathcal{\hat{K}}_{\Sigma}M$. This is
particularly interesting when $M$ is a Calabi-Yau threefold.

There are many interesting questions arise from studying the geometry of
$\mathbb{C}$-VCPs. For example we would like to interpret the
Strominger-Yau-Zaslow mirror transformation for Calabi-Yau manifolds $M$ as
the twistor rotation for the holomorphic symplectic manifolds $\mathcal{\hat
{K}}_{\Sigma}M.$

\bigskip

\textit{Acknowledgments: This paper is partially supported by NSF/DMS-0103355.
Authors express their gratitude to S.L. Kong, X.W. Wang for useful discussions.}

\bigskip

Address: School of Mathematics, University of Minnesota, Minneapolis, MN
55454, USA.

Email: JHLEE@MATH.UMN.EDU, LEUNG@MATH.UMN.EDU
\end{document}